\documentclass[12pt]{article}
\usepackage{color}

\usepackage[italian, english]{babel}
\usepackage[OT2, T1]{fontenc}

\usepackage[utf8]{inputenc}

\usepackage{amssymb,amsmath} 

\allowdisplaybreaks 

\usepackage{hyperref}
\hypersetup{hidelinks}

\newtheorem{theorem}{Theorem}[section]
\newtheorem{proposition}[theorem]{Proposition}
\newtheorem{corollary}[theorem]{Corollary}
\newtheorem{lemma}[theorem]{Lemma}
\newtheorem{remark}[theorem]{Remark}
\newtheorem{definition}[theorem]{Definition}

\def\claim#1.{\noindent {\bf #1.}}

\def\bull{\vrule height 1.8ex width 1.0ex depth .1ex }
\def\flushright#1{{\unskip\nobreak\hfil\penalty50\hskip2em\hbox{}\nobreak\hfil #1\parfillskip=0pt\finalhyphendemerits=0\par}}
\def\QED{\ifmmode\eqno\hbox{$\bull$}\else\flushright{\hbox{$\bull$}}\fi}

\newcommand{\abs}[1]{\left| #1 \right|} 
\newcommand{\norm}[1]{\left\Vert#1\right\Vert} 
\newcommand{\tnorm}[1]{{\left\vert\kern-0.25ex\left\vert\kern-0.25ex\left\vert #1 
		\right\vert\kern-0.25ex\right\vert\kern-0.25ex\right\vert}}
		
\newcommand{\mc}[1]{\mathcal{#1}}

\newcommand{\parag}[1]{\left\{ \begin{aligned} #1 \end{aligned}\right.}

\def\supp{\mathop{\rm supp}}

\newcommand{\Tail}{{\rm Tail}}
\newcommand{\DG}{{\rm DG}}

\def\cat{\mathop{\rm cat}\nolimits}
\newcommand{\cupl}{{\rm cupl}}

\newcommand{\R}{\mathbb{R}}
\newcommand{\N}{\mathbb{N}}
\newcommand{\F}{\mathbb{F}}

\newcommand{\eps}{\varepsilon}

\newcommand{\wto}{\rightharpoonup} 
\newcommand{\nequiv}{\not \equiv}

\usepackage{scalerel}
\usepackage{tikz}
\usetikzlibrary{svg.path}
\definecolor{orcidlogocol}{HTML}{A6CE39}
\tikzset{
 orcidlogo/.pic={
 \fill[orcidlogocol] svg{M256,128c0,70.7-57.3,128-128,128C57.3,256,0,198.7,0,128C0,57.3,57.3,0,128,0C198.7,0,256,57.3,256,128z};
 \fill[white] svg{M86.3,186.2H70.9V79.1h15.4v48.4V186.2z}
 svg{M108.9,79.1h41.6c39.6,0,57,28.3,57,53.6c0,27.5-21.5,53.6-56.8,53.6h-41.8V79.1z M124.3,172.4h24.5c34.9,0,42.9-26.5,42.9-39.7c0-21.5-13.7-39.7-43.7-39.7h-23.7V172.4z}
 svg{M88.7,56.8c0,5.5-4.5,10.1-10.1,10.1c-5.6,0-10.1-4.6-10.1-10.1c0-5.6,4.5-10.1,10.1-10.1C84.2,46.7,88.7,51.3,88.7,56.8z};
 }
}
\newcommand\orcidicon[1]{\href{https://orcid.org/#1}{\mbox{\scalerel*{
\begin{tikzpicture}[yscale=-1,transform shape]
\pic{orcidlogo};
\end{tikzpicture}
}{|}}}}

\begin{document}

\title{On the fractional NLS equation and \\ the effects of the potential well's topology}

\author{
		Silvia Cingolani \orcidicon{0000-0002-3680-9106}\\
		Dipartimento di Matematica\\ 
		Universit\`{a} degli Studi di Bari Aldo Moro
		\\
		Via E. Orabona 4, 70125 Bari, Italy
		\\ \href{mailto:silvia.cingolani@uniba.it}{silvia.cingolani@uniba.it}
		\\ 
		\\
		Marco Gallo \orcidicon{0000-0002-3141-9598}
		\\
		Dipartimento di Matematica\\
		Universit\`{a} degli Studi di Bari Aldo Moro
		\\
		Via E. Orabona 4, 70125 Bari, Italy
		\\ \href{mailto:marco.gallo@uniba.it}{marco.gallo@uniba.it} 
		\\ 
}

\maketitle

\abstract{
\noindent In this paper we consider the fractional nonlinear Schr\"odinger equation
$$ \eps^{2s}(- \Delta)^s v+ V(x) v= f(v), \quad x \in \R^N$$
where $s \in (0,1)$, $N \geq 2$, $V \in C(\R^N,\R)$ is a positive potential and $f$ is a nonlinearity satisfying Berestycki-Lions type conditions. 
For $\eps>0$ small, we prove the existence of at least $\cupl(K)+1$ positive solutions, where $K$ is a set of local minima in a bounded potential well and $\cupl(K)$ denotes the cup-length of $K$. By means of a variational approach, we analyze the topological difference between two levels of an indefinite functional in a neighborhood of expected solutions.
Since the nonlocality comes in the decomposition of the space directly, we introduce a new fractional center of mass, via a suitable seminorm.
Some other delicate aspects arise strictly related to the presence of the nonlocal operator. By using regularity results based on fractional De Giorgi classes, we show that the found solutions decay polynomially and concentrate around some point of $K$ for $\eps$ small.
}

\bigskip
	

\medskip


\newpage

\tableofcontents


\section{Introduction}\label{section:1}
In \cite{Las} Laskin developed a new extension of the fractality concept in quantum physics and formulated the fractional nonlinear Schr\"odinger (fNLS for short) equation 
	\begin{equation}\label{fNLS}
	i \hbar \partial_t \psi = \hbar^{2s} (-\Delta)^s \psi + V(x) \psi- f(\psi), \quad (t,x) \in (0,+\infty) \times \R^N
	\end{equation}
	where $s \in (0,1)$, $N > 2s$, the symbol $(-\Delta)^{s}\psi=\mc{F}^{-1}(|\xi|^{2s} \mc{F}(\psi))$ denotes the fractional power of the Laplace operator defined via Fourier transform $\mc{F}$ on the spatial variable, $\hbar$ designates the usual Planck constant, $V$ is a real potential and $f$ is a Gauge invariant nonlinearity, i.e. $f(e^{i \theta} \rho) = e^{i \theta} f(\rho)$ for any $\rho, \theta \in \R$. 
The nonlocal equation $(\ref{fNLS})$ was derived as extension of the classical NLS equation $(s=1)$, replacing
the path integral over Brownian motions (random motions seen in swirling gas molecules) to L\'evy
 flights (a mix of long trajectories and short, random walk models found in turbulent fluids). 
 The wave function $\psi(x,t)$ represents the quantum mechanical probability amplitude for a given unit mass particle to have position $x$ and time $t$, under the confinement due to the potential $V,$ and $|\psi|^2$ is the corresponding probability density. We refer to \cite{Las, Las3}
 for a detailed discussion of the physical motivation of the {fNLS} equation.
		
Special solutions of the equation $(\ref{fNLS})$ are given by the standing waves, i.e. factorized functions $\psi (t,x) = e^{\frac{i E t}{\hbar}} v(x)$, $E\in \R$. 
For small $\hbar >0$, these standing waves are usually referred to as 
{\sl semiclassical} states and the transition from quantum physics to classical physics is somehow described letting $\hbar \to 0$ (see \cite{SS}).

Without loss of generality, shifting $E$ to
$0$ and denoting $\hbar \equiv \varepsilon$, the search of semiclassical states leads to investigate the following nonlocal equation
\begin{equation} \label{eq:1.1}
\varepsilon^{2s} (-\Delta)^s v + V(x)v = f(v), \quad x \in \R^N 
\end{equation}
when $V$ is positive and $\varepsilon > 0$ is small.

In the limiting case $s=1$ the semiclassical analysis of NLS equations
has been largely investigated, starting from the seminal paper \cite{FW}.
By means of a finite dimensional reduction, Floer and Weinstein proved 	
the existence of positive spike solutions to the 3D cubic NLS equation, 
concentrating at each nondegenerate critical point of the potential $V$ (see also \cite{Oh2}). 
Successively, refined variational techniques were implemented to
study singularly perturbed elliptic problems in entire space: several
 existence results of positive spike solutions to the NLS equation in a semiclassical regime are derived under different assumptions on the potential and the nonlinear terms. 
 We confine to mention 
 \cite{Rab, Wan, DF,ABC,DF2, BJ,BJ2,BT1, BT2}
and references therein.
In \cite{CL,AMS,CJT} topological invariants were used to derive multiplicity results, 
in the spirit of well-known results of Bahri, Coron \cite{BC} and Benci, Cerami \cite{BC1} for semilinear elliptic problems with Dirichlet boundary condition. Precisely, in \cite{CL} it has been proved that the number of positive solutions of the stationary NLS equation is influenced by the topological richness of the set of global minima of $V$. Some years later,
using a perturbative approach, Ambrosetti, Malchiodi and Secchi \cite{AMS} obtained a multiplicity results for the NLS equation with power nonlinearity, assuming that 
the set of critical points of $V$ is nondegenerate in the sense of Bott.
 
More recently, 
in \cite{CJT} Jeanjean, Tanaka and the first author improved the result in \cite{CL}, relating the number of semiclassical standing waves solutions to the 
cup-length of $K$, where $K$ is a set of local (possibly degenerate) minima of the potential, under almost optimal assumptions on the nonlinearity
(see also the recent paper 
\cite{CT}
in the context of nonlinear Choquard equations).

When $s \in (0,1)$, the search of semiclassical standing waves for the fNLS equation has been firstly considered by D\'avila, Del Pino and Wei in \cite{DPW} under the assumptions $f(t)= |t|^{p-2}t$, with $2 < p < 2^*_s$, where $2^*_s: = \frac{2N}{N-2s}$ is the Sobolev critical exponent, and $V \in C^{1,\alpha}(\R^N)$ is bounded. 
Using a Lyapunov-Schmidt reduction inspired by \cite{FW, Oh2}, they showed the existence of a positive spike solution whose maximum point concentrates at some nondegenerate critical point of $V$; this approach relies on the nondegeneracy property of the linearization at the positive ground state shown by Frank, Lenzmann and Silvestre \cite{FLS}.
Successively, inspired by \cite{DF, BJ}, variational techniques were employed to derive existence of spike solutions concentrating at local minima of $V$, 
see \cite{Seo, AlMi, Amb} and references therein (see also \cite{Sec} where global assumptions on $V$ are considered). 
We further cite \cite{FMV} in which necessary conditions are given to the existence of concentrating solutions.

A first multiplicity result for the (fNLS) equation is obtained 
in \cite{FS}, inspired by \cite{CL}. Precisely, letting $K$ be the set of global minima of $V$, Figueiredo and Siciliano proved that the number of positive solutions of $(\ref{eq:1.1})$, when $f$ satisfies monotonicity and Ambrosetti-Rabinowitz condition, is at least given by the Ljusternik-Schnirelmann category of $K$: here the search of solutions of $(\ref{eq:1.1})$ can be reduced to the study of the (global) level sets of the Nehari manifold, where the energy functional is restricted, and to deformation arguments valid on Hilbert manifolds without boundary. 
See also \cite{AA} where the Ambrosetti-Rabinowitz condition is dropped. 
In \cite{Che}, moreover, Chen implemented a Lyapunov-Schmidt reduction for nondegenerate critical points of $V$ and power-type functions $f$ in order to get multiplicity results related to the cup-length, extending the results of \cite{AMS}.

In the present paper we are interested to prove multiplicity of positive solutions for the fNLS equation \eqref{eq:1.1} when $\eps$ is small, without monotonocity and Ambrosetti-Rabinowitz conditions on $f$, nor nondegeneracy and global conditions on $V$.

On the potential $V$ we assume

\begin{itemize}
	\item[(V1)] $\, V\in C(\R^N, \R)\cap L^{\infty}(\R^N)$, $\underline{V}:=\inf_{\R^N} V>0$ (see also Remark \ref{rem_ipotesi_V});
	\item[(V2)] $\,$there exists a bounded domain $\Omega\subset \R^N$ such that
	$$m_0:= \inf_{\Omega} V < \inf_{\partial \Omega} V;$$
\end{itemize}
\vspace{-\topsep}
by the strict inequality and the continuity of $V$, we can assume that $\partial \Omega$ is regular.
We define $K$ as the set of local minima
\begin{equation}\label{eq_def_K}
K:=\{ x \in \Omega \mid V(x)=m_0\}.
\end{equation}

On $f$ we assume
\begin{itemize}
	\item[(f1)] Berestycki-Lions type assumptions \cite{BL} with respect to $m_0$, that is
	\begin{itemize}
		\item[(f1.1)] $\, f\in C(\R, \R)$;
		\item[(f1.2)] $\, \lim_{t \to 0} \frac{f(t)}{t}=0$;
		\item[(f1.3)] $\, \lim_{t \to +\infty} \frac{f(t)}{|t|^p}=0$ for some $p \in (1, 2^*_s-1)$, where we recall $2^*_s=\frac{2N}{N-2s}$; 
		\item[(f1.4)] $\, F(t_0)> \frac{1}{2} m_0 t_0^2$ for some $t_0>0$, where $F(t):= \int_0^t f(s) ds$;
	\end{itemize}

	\item[(f2)] $f(t)=0$ for $t\leq 0$.
\end{itemize}

On $f$ we further assume 
\begin{itemize}
	\item[(f3)] $f\in C^{0,\gamma}_{loc}(\R)$ for some $\gamma \in (1-2s, 1)$ if $s\in (0,1/2]$.
\end{itemize}
We remark that (f3) is needed only to get a Pohozaev identity (see \cite[Proposition 1.1]{BKS}).

\smallskip

Setting $u:= v(\eps \cdot)$, $(\ref{eq:1.1})$ can be rewritten as
\begin{equation}\label{eq_cambio_var}
(-\Delta)^s u + V(\eps x) u= f(u), \quad \textnormal{ $x \in \R^N$},
\end{equation}
thus the equation
\begin{equation}\label{eq_intr_limit_eq}
(-\Delta)^s U + a U = f(U), \quad x \in \R^N,
\end{equation}
for some $a>0$, becomes a formal limiting equation, as $\eps \to 0$, for $(\ref{eq_cambio_var})$.

It is standard that weak solutions to \eqref{eq_cambio_var} correspond to critical points of the $C^1$-energy functional
$$I_{\eps}(u):=\frac{1}{2} \int_{\R^N} \abs{(-\Delta)^{s/2} u}^2 dx+ \frac{1}{2} \int_{\R^N} V(\eps x) u^2 dx- \int_{\R^N}F(u) dx, \quad u\in H^s(\R^N) $$
where $H^s(\R^N)$ is the fractional Sobolev space. 
We remark that, because of the general assumptions on $f$, we can not take advantage of the boundedness of the functional from above and below, nor of Nehari type constraint. 
Therefore in the present paper we combine reduction methods and penalization arguments in a nonlocal setting: 
in particular, as in \cite{CJT, CT}, the analysis of the topological changes between two level sets of the indefinite energy functional $I_{\eps}$ in a small neighborhood $\mc{X}_{\eps, \delta}$
of expected solutions is essential in our approach. 
With the aid of $\eps$-independent pseudo-differential estimates, we detect 
such a neighborhood, which will be positively invariant
under a pseudo-gradient flow, and we develop our deformation argument in the context of nonlocal operators. 
To this aim we introduce two maps $\Phi_\varepsilon$ and $\Psi_\varepsilon$ between topological pairs: we emphasize that to define such maps, a center of mass $\Upsilon$ and a functional $P_a$ which is inspired by the Pohozaev identity are crucial.

With respect to the local case, several difficulties arise linked to special features of the nonlocal nature of the problem: among them we have 
the polynomial decay of the least energy solutions of the limiting problems, the weak regularizing effect of the fractional Laplacian, the lack of general comparison arguments, the differences between the supports of a function and of its Fourier transform, and the lack of the standard Leibniz formula (see e.g. \cite{Sec, BWZ}). 
 
We highlight that, for fractional equations, the nonlocal part strongly influences the decomposition of the space and this makes quite delicate to use truncating test functions and perform the localization of the centers of mass. In the present paper we introduce a new fractional local center of mass by means of suitable seminorm, stronger than the usual Gagliardo seminorm in a bounded set, and we need to implement new ideas to overcome the above obstructions.

\vspace{2mm}
Our main result is the following theorem.
\begin{theorem}\label{claim:1.1}
	Suppose $N\geq 2$ and that \textnormal{(V1)--(V2)}, \textnormal{(f1)--(f3)} hold. 		
	Let $K$ be defined by \eqref{eq_def_K}.
	Then, for sufficiently small $\eps >0$, equation \eqref{eq:1.1} has at least
	$\cupl(K)+1$ positive solutions, which belong to $C^{0, \sigma}(\R^N) \cap L^{\infty}(\R^N)$ for some $\sigma \in (0,1)$.
\end{theorem}

Here $\cupl(K)$ denotes the \emph{cup-length} of $K$ defined by the Alexander-Spanier cohomology with coefficients in some field $\F$ (see Section \ref{sec_cuplenght}). 

\begin{remark}\label{examples}
	Notice that the cup-length of a set $K$ is strictly related to the \emph{category} of $K$. Indeed, if $K=S^{N-1}$ is the $N-1$ dimensional sphere in $\R^N$, then $\cupl(K) + 1= \cat(K) =2$; if $K=T^N$ is the $N$-dimensional torus, then $\cupl(K) + 1 = \cat(K)= N + 1$.
	However in general $\cupl(K) +1 \leq \cat(K)$.
\end{remark}

A multiplicity result in the limiting case $s=1$, similar to Theorem \ref{claim:1.1}, is contained in \cite{CJT}.

\begin{remark}\label{rem_ipotesi_V}
Observe that, arguing as in \cite{BJ,BJ2} and \cite{BT1}, we could omit the assumption that $V$ is bounded from above in Theorem \ref{claim:1.1}. For the sake of simplicity, we assume here the boundedness of $V$.
\end{remark}

By using recent regularity results based on fractional De Giorgi classes and tail functions (see \cite{Coz}) we are able to prove also the following theorem.

\begin{theorem}\label{mean-concen}
	Let $(\eps_n)_{n \in \N}$ with $\eps_n \to 0^+$ as $n \to +\infty$. For sufficiently large $n\in \N$, in the assumptions of Theorem \ref{claim:1.1}, let $v_{\eps_n}$ be one of the $\cupl(K)+1$ solutions of equation \eqref{eq:1.1}.
	Then, up to a subsequence, $(v_{\eps_n})_{n\in \N}$ \emph{concentrates} in $K$ as $n \to +\infty$. More precisely, for each $n\in \N$ there exists a maximum point $x_{\eps_n}\in \R^N$ of $v_{\eps_n}$ such that $$\lim_{n\to +\infty}d(x_{\eps_n}, K) =0.$$
	In addition, $ v_{\eps_n}(\eps_n \cdot+x_{\eps_n})$ converges in $H^s(\R^N)$ and uniformly on compact sets to a least energy solution of
	\begin{equation}\label{eq_least_energy_m0}
	(- \Delta)^s U + m_0 U = f(U), \quad U >0, \quad U \in H^{s}(\R^N),
	\end{equation}
	and, for some positive $C', C''$ independent on $n \in \N$, we have the uniform polynomial decay
	$$\frac{C'}{1+|\frac{x-x_{\eps_n}}{\eps_n}|^{N+2s}}\leq v_{\eps_n}(x) \leq \frac{C''}{1+|\frac{x-x_{\eps_n}}{\eps_n}|^{N+2s}},\quad \textit{ for $x \in \R^N$}.$$
\end{theorem}

\medskip

The paper is organized as follows. In Section \ref{sec_fract_fram} we recall some basic notions on the nonlocal framework. In Section \ref{sez_limit_eq} we overview some results on the limiting equation \eqref{eq_intr_limit_eq} and we introduce
 a new fractional center of mass $\Upsilon$, by means of a suitable seminorm;
 we postpone to Appendix \ref{sez_poly_dec} the proof of the uniform polynomial decay of the solutions of \eqref{eq_intr_limit_eq}, since it shares some argument with the proof of Theorem \ref{mean-concen}. Section \ref{sez_singur_perturb} is the main core of the article, where we introduce a penalized functional and prove a deformation lemma on a neighborhood of expected solutions; moreover, we build suitable maps $\Phi_{\eps}$, $\Psi_{\eps}$ essential in the proof of the multiplicity of solutions. In Section \ref{sec_cuplenght} we prove Theorem \ref{claim:1.1} by the use of the deformation lemma and the built maps applied to the theory of relative category and relative cup-length, of which we briefly recall the definitions. Finally we prove Theorem \ref{mean-concen} by using regularity results based on fractional De Giorgi classes.

\section{Preliminaries}\label{sec_fract_fram}

Let $N>2s$. Set the Gagliardo seminorm
$$[u]_{\R^N}^2:= \int_{\R^N} \int_{\R^N} \frac{|u(x)-u(y)|^2}{|x-y|^{N+2s}} dx \, dy$$
we can define the fractional Sobolev space \cite[Section 2]{DGV}
$$H^s(\R^N):=\left \{ u \in L^2(\R^N) \mid [u]_{\R^N} < + \infty\right\}$$
and set for every $u \in H^s(\R^N)$
$$\norm{u}_{H^s(\R^N)}^2 := \norm{u}_{L^2(\R^N)}^2 + [u]_{\R^N}^2;$$
here $\norm{\cdot}_{L^q(A)}$ denotes the $L^q$-Lebesgue norm in $A\subset \R^N$ for $q \in [1, +\infty]$, and we will often write $\norm{\cdot}_q:=\norm{\cdot}_{L^q(\R^N)}$.
Moreover, we define the fractional Laplacian
\begin{equation}\label{eq_def_frac_lapl}
(-\Delta)^s u(x):= C \int_{\R^N} \frac{u(x)-u(y)}{|x-y|^{N+2s}} dx, \quad u \in H^s(\R^N),
\end{equation}
where $C$ (as the following appearing constants) depends on $N$ and $s$, and the integral is in the principal value sense; the following equality holds \cite[Proposition 3.6]{DGV}
\begin{equation}\label{eq_equiv_norm_laplac}
\norm{(-\Delta)^{s/2}u}_2^2 = C' [u]_{\R^N}^2 
\end{equation}
and more generally (by polarization identity)
$$ \int_{\R^N} (-\Delta)^{s/2} u \, (-\Delta)^{s/2} v \, dx = C'\int_{\R^N} \int_{\R^N} \frac{\big(u(x)-u(y)\big)\big(v(x)-v(y)\big)}{|x-y|^{N+2s}} dx \, dy .$$
We further have \cite[Theorem 6.5]{DGV} 
$$H^s(\R^N) \hookrightarrow L^q(\R^N)$$
for $q \in [2, 2^*_s]$ 
and, for $u \in H^s(\R^N)$,
\begin{equation}\label{eq_emb_homog}
\norm{u}_{L^{2^*_s}(\R^N)}\leq C'' \norm{(-\Delta)^{s/2}u}_{L^2(\R^N)}.
\end{equation}
Moreover, if $q\in [2,2^*_s)$ we have \cite[Corollary 7.2]{DGV}
$$H^s(\R^N) \hookrightarrow \hookrightarrow L^q_{loc}(\R^N)$$
in the sense that for every $(u_n)_n$ bounded in $H^s(\R^N)$, and for every $A\subset \R^N$ bounded and regular enough ($\partial A$ Lipschitz), we have that $(u_n)_n$ restricted to $A$ admits a convergent subsequence in $L^q(A)$. 

We will make use also of
$$\norm{u}_{H_{\eps}^s(\R^N)}^2:= \norm{(-\Delta)^{s/2} u}_2^2 + \int_{\R^N} V(\eps x) u^2 dx$$
which is an equivalent norm on $H^s(\R^N)$, thanks to the positivity and the boundedness of $V$ and \eqref{eq_equiv_norm_laplac}; the space $H^s_{\eps}(\R^N)$ is defined straightforwardly.

Finally, we will make use of the following \emph{mixed} Gagliardo seminorm 
$$[u]_{A_1, A_2}^2:= \int_{A_1} \int_{A_2} \frac{|u(x)-u(y)|^2}{|x-y|^{N+2s}} dx \, dy, \quad [u]_A:=[u]_{A,A}$$
for any $A_1, A_2,A \subset \R^N$ and $u \in H^s(\R^N)$;
by using that $\varphi_u(x,y):= \frac{|u(x)-u(y)|}{|x-y|^{N/2+s}}$ satisfies $\varphi_{u+v} \leq \varphi_u + \varphi_v$ and $[u]_{A_1, A_2} = \norm{\varphi_u}_{L^2(A_1 \times A_2)}$, we have that $[u]_{A_1, A_2}$ is actually a seminorm. 

For any $u \in H^s(\R^N)$ and $A \subset \R^N$ it will be useful to work also with the following norms:
\begin{equation}\label{eq_def_Htilde}
\norm{u}_{A}^2:= \norm{u}_{L^2(A)}^2 + [u]_{A,\R^N}^2
\end{equation}
and
\begin{equation}\label{eq_def_triplanorm}
\tnorm{u}_{A}:=\norm{u}_{L^{p+1}(A)} + \norm{u}_{A},
\end{equation}
where $p$ is introduced in assumption (f1.3). We highlight that $\norm{u}_{\R^N}=\norm{u}_{H^s(\R^N)}$, but generally $\norm{u}_A \geq \norm{u}_{H^s(A)}$ for $A \neq \R^N$ (see \cite[Section 2]{DGV} for a definition of $H^s(A)$). 
By $H^s(A) \hookrightarrow L^{p+1}(A)$ \cite[Theorem 6.7]{DGV} the norms $\norm{\cdot}_A$ and $\tnorm{\cdot}_A$ are equivalent: on the other hand, the constant such that $\tnorm{u}_A \leq C_A \norm{u}_A$ depends on $A$, thus not useful for $\eps$-dependent sets $A=A(\eps)$. This is why we will make direct use of $\tnorm{\cdot}_A$. 

Moreover, we denote by 
$$[u]_{C^{0,\sigma}(A)}:= \sup_{\substack{x, y \in A \\ x\neq y}} \frac{\abs{u(x)-u(y)}}{\abs{x-y}^{\sigma}}$$
the usual seminorm in H\"older spaces for $\sigma \in (0,1]$. We allow $\sigma>1$ by simply writing $C^{\sigma}(\R^N)$.

Before ending this section, we highlight that the assumptions on the function $f$ imply the following standard fact, which will be extensively used throughout the paper: 
for each $q\geq p$ and $\beta>0$ there exists a $C_{\beta}>0$ such that
\begin{equation}\label{eq_prop_f}
\abs{f(t)}\leq \beta |t| + C_{\beta} |t|^q \quad \textnormal{ and } \quad \abs{F(t)}\leq C\left(\beta |t|^2 + C_{\beta} |t|^{q+1}\right).\end{equation}

\medskip

\textbf{Further notation}. We highlight that, all throughout the paper, we will assume $N\geq 2$, and the constants $C, C'$ appearing in inequalities may change from a passage to another. 
To avoid cumbersome notations, we will not stress the dependence of such constants, which will be based only on the fixed quantities. 
Moreover, by $\approx$ we will mean \emph{approximately equal to}.

We will write $B_R(x_0)$ for the ball, centered in $x_0 \in \R^N$, with radius $R>0$; in particular, $B_R:=B_R(0)$. Moreover 
$$A_{\delta}:= \{x \in X \mid d(x,A)\leq \delta\}$$
for any $A\subset (X,d)$ metric space. In addition, we will sometimes write $\complement(A):= \R^N \setminus A$ for $A \subset \R^N$ to avoid cumbersome notation. 

Finally, for every $A\subset B \subset \R^N$, we will write
$$A\prec \phi \prec B$$
to indicate a Urysohn-type regular function $\phi \in C^{\infty}_c(\R^N)$ such that
$$\phi_{|A}=1 \quad \textnormal{ and } \quad \phi_{|\R^N \setminus B}=0.$$

\section{Limiting equation}\label{sez_limit_eq}

\subsection{A single equation} \label{subsec_single_eq}

Consider
\begin{equation}\label{eq_limit_eq}
(-\Delta)^s U + a U = f(U), \quad x \in \R^N
\end{equation}
with $a>0$. Weak solutions of \eqref{eq_limit_eq} are known to be characterized as critical points of the $C^1$-functional $L_a: H^s(\R^N) \to \R$
$$L_a(U):= \frac{1}{2} \norm{(-\Delta)^{s/2}U}_2^2 + \frac{a}{2} \norm{U}_2^2 - \int_{\R^N} F(U) dx, \quad U \in H^s(\R^N).$$
Set moreover the Pohozaev functional $P_a: H^s(\R^N) \setminus \{0\} \to \R_+$
$$P_a(U):= \left( \frac{2N}{N-2s} \; \frac{\int_{\R^N} F(U) dx - \frac{a}{2} \norm{U}_2^2}{\norm{(-\Delta)^{s/2} U}_2^2}\right)_+^{\frac{1}{2s}}, \quad U\in H^s(\R^N), \; U \neq 0.$$
By \cite[Proposition 1.1]{BKS} we have that the Pohozaev identity
$$P_a(U)=1$$
holds for each nontrivial solution $U$ of \eqref{eq_limit_eq}. To reach this claim we need the solutions to be regular enough, fact that is given by (f3). 
We observe that the functional $P_a$ will be of key importance for estimating $L_a$ from below, see Lemma \ref{lemma_stima_g0} and Lemma \ref{lemma_stima_g}. 

In addition, we recall the following result by \cite[Theorem 1.2]{BKS}.
\begin{theorem}[Existence of a Pohozaev minimum]\label{thm_esist_Poh_min}
Assume \textnormal{(f1)} with respect to $a>0$ and \textnormal{(f2)}. Let
$$C_{po,a}:=\inf \big\{ L_a(U) \mid U \in H^s(\R^N) \setminus \{0\}, \; P_a(U)=1\big\}.$$
Then there exists a positive minimizer, which is a weak solution of \eqref{eq_limit_eq}; this minimizer is generally not unique, even up to translation. Moreover, $C_{po,a}>0$.
\end{theorem}

We further denote by
$$E_a:=\inf \big\{L_a(U) \mid U \in H^s(\R^N) \setminus \{0\}, \;  L'_a(U)=0 \big\}$$
the least energy for $L_a$. As an easy consequence of Theorem \ref{thm_esist_Poh_min} we have
\begin{equation}\label{eq_ug_poh_gr.st}
E_a = C_{po,a}.
\end{equation}

In the following result we recall \cite[Theorem 1.1 and Theorem 1.3]{BKS} and \cite[Theorem 1.5]{FQT}, which state the regularity and the polynomial decay of solutions of \eqref{eq_limit_eq}. This decay is much less slower than the one, exponential, of the local case $s=1$. An alternative proof, which underlines some uniformity in $a$, can be found in the Appendix (see Proposition \ref{prop_polyn_dec}).

\begin{proposition}[Regularity and polynomial decay]\label{thm_reg_dec}
Assume \textnormal{(f1)} with respect to $a>0$ and \textnormal{(f2)}. Then every weak solution $U \in H^s(\R^N)$ of \eqref{eq_limit_eq} is actually a strong solution, i.e. $U$ satisfies \eqref{eq_limit_eq} almost everywhere. Moreover
$$U\in H^{2s}(\R^N)\cap C^{\sigma}(\R^N)$$
for every $\sigma \in (0,2s)$, and
$$U>0 \quad \textnormal{and} \quad \lim_{|x|\to +\infty} U(x)=0.$$
Furthermore, we have a polynomial decay at infinity, that is there exist positive constants $C_a', C_a''$ such that
\begin{equation}\label{eq_bound_funct}
\frac{ C_a'}{1+|x|^{N+2s}} \leq U(x)\leq \frac{ C_a''}{1+|x|^{N+2s}},\quad \textit{ for $x \in \R^N$}.
\end{equation}
\end{proposition}

We observe that the bounding functions in \eqref{eq_bound_funct} belong to $L^q(\R^N)$ for any $q\in [1,+\infty]$. 

\smallskip

We end this section by a technical lemma, which allows to link the level $L_{a}(u)$ of a whatever function $u$, having a functional $P_a(u) \approx 1$, with the ground state $E_{a}$; in particular, it provides a useful lower bound for $L_a$.

\begin{lemma}\label{lemma_stima_g0}
Let $u\in H^s(\R^N)$ and define
$$g(t):=\frac{1}{2s} \left(N t^{N-2s} - (N-2s) t^N\right), \quad t \in \R.$$
\begin{itemize}
\item[(a)] If $P_a(u)\in \Big(0, \left(\frac{N}{N-2s}\right)^{\frac{1}{2s}}\Big)$, which we highlight is a neighborhood of $1$, then
$$L_a(u)\geq g(P_a(u))E_a.$$
\item[(b)] If $u=U_a\left(\frac{\cdot-q}{t}\right)$ for some $q \in \R^N$ and $t \in \R$, with $U_a$ being a ground state of \eqref{eq_limit_eq}, then the above inequality is indeed an equality, that is
$$ L_a\left ( U_a\left(\frac{\cdot-q}{t}\right) \right) = g(t)E_a. $$
\end{itemize}
We highlight that the function $g$ verifies
$$g(t)\leq 1 \quad \textnormal{ and } \quad g(t)=1 \iff t=1.$$
\end{lemma}

\claim Proof.
Let $\sigma:=P_a(u)$ and set $v:=u(\sigma \cdot)$. Then $P_a(v)=1$. A straightforward computation shows, by using that $P_a(v)=1$ and that $g(\sigma)>0$,
$$L_a(u) = g(\sigma) L_a(v) \geq g(\sigma) C_{po,a} = g(\sigma) E_a.$$
We see that, if $u=U_a\left(\frac{\cdot-q}{t}\right)$, then by $P_a(U_a)=1$ we have $\sigma=P_a\left( U_a\left(\frac{\cdot-q}{t}\right) \right)= t$ and thus $v=U_a(\cdot-q)$, which by translation invariance is again a ground state of \eqref{eq_limit_eq}; thus $L_a(v)=C_{po,a}$. This concludes the proof.
\QED

\subsection{A family of equations}
\label{subsec_family}

In this section we study equation \eqref{eq_limit_eq} for variable values of $a>0$. Introduce the notation
$$\Omega[a,b]:=V^{-1}([m_0+a,m_0+b])\cap \Omega$$
and similarly $\Omega(a,b)$ and mixed-brackets combinations. 
We choose now a small $\nu_0>0$ such that the minimum $m_0$ is not heavily perturbed, namely
\begin{itemize}
\item Berestycki-Lions type assumptions (f1) hold with respect to $a \in [m_0, m_0+\nu_0]$ -- i.e., in particular, $F(t_0)> \frac{1}{2} (m_0+\nu_0) t_0^2$;
\item assumption (V2) holds with respect to $m_0+\nu_0$, i.e. $m_0 + \nu_0 < \inf_{\partial \Omega} V$;
\item $\overline{\Omega[0,\nu_0]}\subset K_d \subset \Omega$ for a sufficiently small $d>0$ subsequently fixed, see Lemma \ref{lem_cupl_jK};
\item other conditions subsequently stressed, see e.g. \eqref{eq_stima_l0} and Lemma \ref{lemma_stima_g}.
\end{itemize}
We observe that, by construction, for $a \in [m_0, m_0+\nu_0]$ the considerations of Section \ref{subsec_single_eq} apply. Moreover, by scaling arguments on $C_{po,a}$, we notice that $a\in [m_0, m_0+\nu_0] \mapsto E_a \in (0, +\infty)$ is strictly increasing and that, up to choosing a smaller $\nu_0$, we have $E_{m_0+\nu_0} < 2E_{m_0}$ and thus we can find an $l_0=l_0(\nu_0)\in \R$ such that
\begin{equation}\label{eq_stima_l0}
E_{m_0+\nu_0} < l_0 < 2E_{m_0}.
\end{equation}
At the end of the paper, we will make $\nu_0$ and $l_0$ moving such that $\nu_0^n \to 0$ and $l_0^n\to E_{m_0}$. 
We now define the set of \emph{almost} ground states of \eqref{eq_limit_eq}
$$S_a:= \left \{ U \in H^s(\R^N)\setminus\{0\} \mid L'_a(U)=0,\; L_a(U)\leq l_0, \; U(0)=\max_{\R^N} U\right \}\neq \emptyset.$$
We observe that we set the last condition in order to fix solutions in a point and prevent them to escape to infinity; the idea is to gain thickness and compactness (see \cite{BJ, CSS, CJT}). 
We further define
$$\widehat{S}:= \bigcup_{a\in [m_0, m_0+\nu_0]} S_a.$$
The following properties of the set $\widehat{S}$ will be of key importance in the whole paper.

\begin{lemma} \label{lemma_stima_unif} 
The following properties hold.
\begin{itemize}
\item[(a)] There exist positive constants $C', C''$ such that, for each $U\in \widehat{S}$ we have
\begin{equation}\label{eq_stima_uniformeU}
\frac{C'}{1+|x|^{N+2s}} \leq U(x)\leq \frac{C''}{1+|x|^{N+2s}},\quad  \textit{ for $x \in \R^N$}.
\end{equation}
\item[(b)] $\widehat{S}$ is compact. Since it does not contain the zero function, we have
$$r^*:= \min_{U \in \widehat{S}} \norm{U}_{H^s(\R^N)} >0;$$
the maximum is attained as well.
\item[(c)] We have
$$\lim_{R\to +\infty} \norm{U}_{\R^N \setminus B_R} =0, \quad \textit{uniformly for $U\in \widehat{S}$},$$
where the norm $\norm{\cdot}_{\R^N\setminus B_R}$ is defined in \eqref{eq_def_Htilde}. 
Moreover, if $(U_n)_n \subset \widehat{S}$, and $(\theta_n)_n \subset \R^N$ is included in a compact set, then
$$\lim_{n \to +\infty} \norm{U_n(\cdot + \theta_n)}_{\R^N \setminus B_n}=0.$$
\end{itemize}
\end{lemma}

\claim Proof.

\textbf{Step 1.} We see that $\widehat{S}$ is bounded. Indeed, by the Pohozaev identity, we have
$$\norm{(-\Delta)^{s/2}U}_2^2 = \frac{N}{s} L_a(U) \leq \frac{N}{s} l_0.$$
By \eqref{eq_emb_homog} we have that also $\norm{U}_{2^*_s}$ is uniformly bounded. 
Since $L'_a(U)U=0$, we have by \eqref{eq_prop_f}
$$\norm{(-\Delta)^{s/2}U}_2^2 + a\norm{U}_2^2 = \int_{\R^N} f(U)U dx\leq \beta \norm{U}_2^2 + C_{\beta} \norm{U}_{2^*_s}^{2^*_s}$$
which implies, by chosing $\beta < a$, that also $\norm{U}_2^2$ is bounded.

\textbf{Step 2.} There exist uniform $C>0$ and $\sigma \in (0,1)$ such that
\begin{equation}\label{eq_dim_stime_uniformi}
\norm{U}_{\infty} \leq C, \quad [U]_{C^{0,\sigma}_{loc}(\R^N)} \leq C
\end{equation}
for any $U \in \widehat{S}$. We postpone to the Appendix the proof of \eqref{eq_dim_stime_uniformi}, as well as the proof of the uniform pointwise estimate \eqref{eq_stima_uniformeU}; see Proposition \ref{prop_polyn_dec}.

\medskip

We show now $(b)$, which is an adaptation of the fact that $S_a$ itself is compact.

\textbf{Step 3.} We observe first that $\widehat{S}$ is closed. Indeed, if $U_k \in S_{a_k} \subset \widehat{S}$ converges strongly to $U$, then up to a subsequence we have $a_k \to a\in [m_0,m_0+\nu_0]$ and, by the strong convergence, we have that the condition
$$E_{m_0} \leq L_a(U)\leq l_0 $$
holds, which in particular implies that $U\nequiv 0$.
Moreover, exploiting the weak convergence $U_k \wto U$, we obtain by classical arguments that for each $v\in H^s(\R^N)$
\begin{eqnarray*}
\lefteqn{0=L'_{a_k}(U_k)v = \int_{\R^N} (-\Delta)^{s/2} U_k (-\Delta)^{s/2} v \, dx + a_k \int_{\R^N} U_k v \, dx - \int_{\R^N} f(U_k) v \, dx} \\
&\to& \int_{\R^N} (-\Delta)^{s/2} U(-\Delta)^{s/2} v \, dx+ a \int_{\R^N} U v \, dx- \int_{\R^N} f(U) v \, dx= L'_a(U)v,
\end{eqnarray*}
that is, $L'_a(U)=0$. As regards the maximum in zero, we need a pointwise convergence. In order to get it, we exploit the fact that, by \eqref{eq_dim_stime_uniformi}, $U_k$ are uniformly bounded in $L^{\infty}(\R^N)$ and in $C^{0, \sigma}_{loc}(\R^N)$ and we apply Ascoli-Arzel\`{a} theorem to get local uniform convergence. 
This shows that $U\in S_a \subset \widehat{S}$.

\textbf{Step 4.} Let now $U_k \in S_{a_k} \subset \widehat{S}$. 
By the boundedness, up to a subsequence we have $U_k \wto U \in H^s(\R^N)$ and $a_k \to a \in [m_0,m_0+\nu_0]$. We need to show that $\norm{U_k}_{H^s(\R^N)}\to \norm{U}_{H^s(\R^N)}$; the closedness of $\widehat{S}$ will conclude the proof.

As observed in Step 3, we have by the weak convergence that $L'_{a_k}(U_k) U_k=0=L'_a(U)U$; hence, if $R>0$ is a radius to be fixed, we gain
\begin{eqnarray*}
\lefteqn{\left |\left(\norm{(-\Delta)^{s/2}U_k}_2^2 + a_k \norm{U_k}_2^2 \right) - \left(\norm{(-\Delta)^{s/2}U}_2^2 + a \norm{U}_2^2 \right) \right|}\\
 &\leq& \Bigg| \int_{|x| \leq R} f(U_k)U_k dx - \int_{|x|\leq R} f(U) U dx \Bigg | + \\
 &&+ \int_{|x| > R} |f(U_k) U_k| dx + \int_{|x|>R} |f(U)U)| dx =: (I) + (II).
\end{eqnarray*}
Fix now a small $\eta>0$. As regards $(II)$, we have by \eqref{eq_prop_f} (recall that the definition of $p$ is given in condition (f1.3))
$$\int_{|x|>R} |f(U_k)U_k| dx\leq \beta \int_{|x|>R} |U_k|^2 dx+ C_{\beta} \int_{|x|>R} |U_k|^{p+1} dx< \eta$$
for sufficiently (uniformly in $k$) large $R>0$ thanks to \eqref{eq_stima_uniformeU}; up to taking a larger $R$, it holds also for $U$.

Fixed this $R>0$, focusing on $(I)$, by classical arguments we have 
$$ \Bigg| \int_{|x| \leq R} f(U_k)U_k dx - \int_{|x|\leq R} f(U) U dx \Bigg |<\eta$$
for sufficiently large $k$. 
Merging together, we obtain
$$\norm{(-\Delta)^{s/2}U_k}_2^2 + a_k \norm{U_k}_2^2 \to \norm{(-\Delta)^{s/2}U}_2^2 + a \norm{U}_2^2$$
which with elementary passages leads to the claim. 

\medskip

\textbf{Step 5.}
Finally, we prove $(c)$. By contradiction, there exists an $\eta>0$ such that, for each $n \in \N$ there exists $U_n \in \widehat{S}$ which satisfies
$$\norm{U_n}_{\R^N \setminus B_n} > \eta.$$
By the compactness, we have, up to a subsequence, $U_n \to U \in \widehat{S}$ as $n \to +\infty$. 
Thus (notice that $\int_{\R^N} \frac{|U(x)-U(\cdot)|^2}{|x-\cdot|^{N+2s}} dx \in L^1(\R^N)$ and absolute integrability of the integral applies)
$$\eta < \norm{U_n}_{\R^N \setminus B_n} \leq \norm{U_n -U}_{H^s(\R^N)} + \norm{U}_{\R^N \setminus B_n} \to 0$$
which is an absurd. 

For the second part, we argue similarly. Indeed, up to a subsequence, $U_n \to U$ in $H^s(\R^N)$ and $\theta_n \to \theta$ in $\R^N$, thus
\begin{eqnarray*}
\lefteqn{\norm{U_n(\cdot-\theta_n)}_{\R^N\setminus B_n}}\\
&\leq& \norm{U_n- U}_{H^s(\R^N)} +\norm{\tau_{\theta_n}U - \tau_{\theta}U}_{H^s(\R^N)} + \norm{U(\cdot-\theta)}_{\R^N\setminus B_n}
\to0,
\end{eqnarray*}
where $\tau_{\theta}$ is the translation. This concludes the proof.
\QED

\medskip

Gained compactness, we turn back considering the set of all the solutions (with no restrictions in zero), that is
$$\widehat{S}':= \bigcup_{p\in \R^N} \tau_p(\widehat{S});$$
we observe that $\widehat{S}'$ is bounded. Moreover we define an open $r$-neighborhood of $\widehat{S}'$, reminiscent of the perturbation approach in \cite{FW, ABC, DPW},
$$S(r):= \{ u \in H^s(\R^N) \mid d(u,\widehat{S}')< r\}$$ 
that is
\begin{eqnarray*}
\lefteqn{S(r)=}\\
&&\Big\{ u=U(\cdot-p)+\varphi \in H^s(\R^N) \mid U \in \widehat{S}, \; p \in \R^N, \; \varphi \in H^s(\R^N), \; \norm{\varphi}_{H^s(\R^N)}< r \Big\}.
\end{eqnarray*}
We further define a \emph{minimal radius} map $\widehat{\rho}: H^s(\R^N) \to \R_+$ by
$$\widehat{\rho}(u):= \inf \left \{ \norm{u-U(\cdot-y)}_{H^s(\R^N)} \mid U \in \widehat{S}, \; y \in \R^N \right\}.$$
We observe
\begin{equation}\label{eq_stima_widerho}
u \in S(r) \implies \widehat{\rho}(u) < r,
\end{equation}
and in addition
$$\widehat{\rho}(u)=\inf\{t\in \R_+ \mid u \in S(t)\},$$
where the infimum on the right-hand side is not attained. Finally, $\widehat{\rho}\in Lip(H^s(\R^N), \R)$ with Lipschitz constant equal to $1$, that is, for every $u,v \in H^s(\R^N)$,
\begin{equation}\label{eq_lipsc_rho}
\abs{\widehat{\rho}(u)-\widehat{\rho}(v)} \leq \norm{u-v}_{H^s(\R^N)}.
\end{equation}

We end this section with two technical lemmas. The first one is a direct consequence of Lemma \ref{lemma_stima_g0}, and allows to link the level $L_{m_0}(u)$ of a whatever function $u\in S(r)$ with the ground state $E_{m_0}$, once $r$ is sufficiently small; this further gives a lower bound for the functional $L_{m_0}$.

\begin{lemma}\label{lemma_stima_g}
Up to taking a smaller $\nu_0=\nu_0(l_0)>0$,
there exists a sufficiently small $r'=r'(\nu_0, r^*)>0$ such that, for every $u \in S(r')$, we have
$$L_{m_0}(u)\geq g(P_{m_0}(u))E_{m_0}.$$
\end{lemma}

\claim Proof.
By Lemma \ref{lemma_stima_g0} (a), we know that the inequality holds if $P_{m_0}(u)$ is in a neighborhood of the value $1$. Observe that $P_{a}(U)=1$ if $U\in S_{a}$: by continuity and compactness, $P_{m_0}(U) \approx 1$ if $U\in S_a$ and $a \approx m_0$. In particular, by choosing a small value of $\nu_0$, $P_{m_0}(U) \approx 1$ for $U \in \widehat{S}$. Indeed
$$
P_{m_0}(U)=1+\frac{1}{2s} \frac{N}{N-2s} \; (a-m_0)\frac{ \norm{U}_2^2}{\norm{(-\Delta)^{s/2} U}_2^2} + o(1);
$$
the addendum on the right-hand side can be bounded by the maximum and the minimum over $\widehat{S}$ (notice that $(-\Delta)^{s/2}U$ cannot be zero) and thus we can find a uniform small $\nu_0$. 
Again by continuity we have that $P_{m_0}(u) \approx 1$ for $u\in S(r')$, $r'$ sufficiently small. Indeed
$$
P_{m_0}(u)= 1-\frac{1}{2s} \frac{N}{N-2s}\frac{1}{\norm{(-\Delta)^{s/2}U}_2^2}\big(L_{m_0}(U(\cdot-p)+\varphi)-L_{m_0}(U(\cdot-p)) \big)+ o(1).
$$
This concludes the proof.
\QED

\bigskip

We notice that the condition $E_{m_0+\nu_0}< l_0$ keeps holding by decreasing $\nu_0$, so no ambiguity in the $l_0$-depending choice of $\nu_0$ in Lemma \ref{lemma_stima_g} arises. 
We focus now on the second lemma. 

\begin{lemma}\label{lem_tec2}
There exist $\nu_1 \in (0, \nu_0)$ and $\delta_0=\delta_0(\nu_1)>0$ such that
$$L_{m_0+\nu_1}(U) \geq E_{m_0}+\delta_0, \quad \textit{uniformly for $U \in \widehat{S}$}.$$
\end{lemma}

\claim Proof.
Observe first that, since $\widehat{S}$ is compact also in $L^2(\R^N)$, we have also finite and strictly positive minimum $\underline{M}$ and maximum $\overline{M}$ with respect to $\norm{\cdot}_2$. Consider $\nu_1 \in (0, \nu_0)$ such that 
$$E_{m_0+\nu_1}-E_{m_0} > \frac{1}{2}(\nu_0-\nu_1)\overline{M};$$
we notice that such $\nu_1$ exists since, as $\nu_1\to \nu_0^+$, the left hand side positively increases while the right hand side goes to zero. Let now $a\in [m_0, m_0+\nu_0]$; we consider two cases. 
If $a\in [m_0,m_0+\nu_1]$ we argue as follow: for $U\in S_a$ we have
\begin{eqnarray*}
L_{m_0+\nu_1}(U) &=& L_a(U) + \tfrac{1}{2}(m_0 + \nu_1 - a)\norm{U}_2^2 \\
&\geq& E_a + \tfrac{1}{2}(m_0 + \nu_1 - a)\underline{M}=:(I)+(II);
\end{eqnarray*}
now, the quantity $(I)$ is minimum when $a=m_0$, while $(II)$ is minimum when $a=m_0+\nu_1$; if both could apply at the same time, we would have as a minimum the quantity $E_{m_0}$. Since it is not possible, we obtain
$$\inf_{U\in \bigcup_{a\in [m_0, m_0+\nu_1]} S_a} L_{m_0+\nu_1}(U) > E_{m_0}.$$
If $a\in (m_0+\nu_1, m_0+\nu_0]$ instead, we have
$$ L_{m_0+\nu_1}(U) \geq E_a - \tfrac{1}{2}(a-(m_0 + \nu_1))\overline{M} \geq E_{m_0+\nu_1} - \tfrac{1}{2}(\nu_0-\nu_1)\overline{M} $$
and thus, by the property on $\nu_1$,
$$\inf_{U\in \bigcup_{a\in [m_0+\nu_1, m_0+\nu_0]} S_a} L_{m_0+\nu_1}(U) \geq E_{m_0+\nu_1} - \tfrac{1}{2}(\nu_0-\nu_1)\overline{M}>E_{m_0}.$$
This concludes the proof, by taking as $\delta_0>0$ the smallest of the two differences.
\QED

\subsection{Fractional center of mass}

As in \cite{CJT}, inspired by \cite{BC2, DF2, BT1}, we want to define a \emph{barycentric map} $\Upsilon$ which, given a function $u=U(\cdot-p) + \varphi \in S(r)$, gives an estimate on the maximum point $p$ of $U(\cdot-p)$; since $\varphi$ is small, $p$ is, in some ways, the \emph{center of mass} of $u$. The idea will be to bound $\Upsilon(u)$ in order to re-gain compactness.

Since the nonlocality comes into the very definition of the ambient space, we need the use of the norm \eqref{eq_def_Htilde}, which we notice being stronger than the one induced by the Gagliardo seminorm.

\begin{lemma} \label{lem_def_bar}
Let $r^*$ be as in Lemma \ref{lemma_stima_unif}. Then there exist a sufficiently large $R_0>0$, a sufficiently small radius $r_0\in (0, r^*)$ and a continuous map
$$\Upsilon: S(r_0) \to \R^N$$
such that, for each $u=U(\cdot-p)+\varphi \in S(r_0)$ we have
$$|\Upsilon(u)-p|\leq 2R_0.$$
Moreover, $\Upsilon$ is continuous and
$-\Upsilon$ is \emph{shift-equivariant}, that is, $\Upsilon(u(\cdot+\xi))=\Upsilon(u)-\xi$ for every $u \in S(r_0)$ and $\xi\in \R^N$.
\end{lemma}

\claim Proof.
Recalled that $r^*=\min_{U\in \widehat{S}} \norm{U}_{H^s(\R^N)}>0$, we have by Lemma \ref{lemma_stima_unif}
\begin{equation}\label{eq_dim_cent_mass}
\norm{U}_{\R^N\setminus B_{R_0}} <\tfrac{1}{8}r^*, \quad \textnormal{ uniformly for $U \in \widehat{S}$}
\end{equation}
for $R_0 \gg 0$. Thus
$$r^* \leq \norm{U}_{H^s(\R^N)} \leq \norm{U}_{B_{R_0}} + \norm{U}_{\R^N\setminus B_{R_0}} < \norm{U}_{B_{R_0}} + \tfrac{1}{8}r^*$$
which implies
$$\norm{U}_{B_{R_0}} > \tfrac{7}{8} r^* \quad \textnormal{and} \quad \norm{U}_{\R^N\setminus B_{R_0}}<\tfrac{1}{8} r^* $$
for each $U \in \widehat{S}$. Consider now a cutoff function $\psi \in C^{\infty}_c(\R_+)$ such that
$$[0,\tfrac{1}{4} r^*] \prec \psi \prec [\tfrac{1}{2} r^*,+\infty);$$
let $r_0\in (0, \frac{1}{8} r^*)$ and define, for each $u\in S(r_0)$ and $q\in \R^N$, a \emph{density} function
{\tiny }$$d(q,u):= \psi \left(\inf_{\tilde{U}\in \widehat{S}} \norm{u-\tilde{U}(\cdot-q)}_{B_{R_0}(q)}\right).$$
Notice that $d(\cdot, u)$ is an integrable function: indeed $[u(\cdot + \xi)]_{B_{R_0}(q)} = [u]_{B_{R_0}(q+\xi)}$ for $\xi \in \R^N$, and $q\mapsto \norm{u-\tilde{U}(\cdot-q)}_{B_{R_0}(q)}=\norm{\tau_q u-\tilde{U}}_{B_{R_0}}$ is continuous by
$$\abs{\norm{\tau_q u-\tilde{U}}_{B_{R_0}}-\norm{\tau_pu -\tilde{U}}_{B_{R_0}}} \leq \norm{\tau_q u - \tau_p u}_{B_{R_0}} \leq \norm{\tau_q u - \tau_p u}_{H^s(\R^N)} \to 0$$
as $p\to q$, and the infimum over continuous functions is upper semicontinuous. 

If we show that $d(\cdot,u)\geq 0$ is not identically zero and it has compact support, then it will be well defined the quantity
$$\Upsilon(u):=\frac{\int_{\R^N} q \, d(q,u) dq}{\int_{\R^N} d(q,u) dq}.$$
We show first that $d(\cdot,u)$ has compact support. Indeed if $u=U(\cdot-p)+\varphi$ and $\tilde{U}\in \widehat{S}$ is arbitrary, then
\begin{eqnarray*}
\lefteqn{\norm{u-\tilde{U}(\cdot-q)}_{B_{R_0}(q)}} \\
&\geq & \norm{\tilde{U}(\cdot-q)}_{B_{R_0}(q)} - \norm{U(\cdot-p)}_{B_{R_0}(q)} - \norm{\varphi}_{B_{R_0}(q)} \\
&\geq & \norm{\tilde{U}}_{B_{R_0}} - \norm{U}_{B_{R_0}(q-p)} -\norm{\varphi}_{H^s(\R^N)} \geq \tfrac{6}{8} r^* - \norm{U}_{B_{R_0}(q-p)};
\end{eqnarray*}
take now $q\notin B_{2R_0}(p)$: if $x \in B_{R_0}(q-p)$, by the fact that $|x-(q-p)|<R_0$ and $|q-p|\geq 2R_0$, we obtain that $\abs{x} \geq R_0$, that is, $x \in \R^N \setminus B_{R_0}$. Therefore by \eqref{eq_dim_cent_mass}
\begin{equation}\label{eq_dim_58} 
\norm{u-\tilde{U}(\cdot-q)}_{B_{R_0}(q)} \geq \tfrac{6}{8} r^* - \norm{U}_{\R^N \setminus B_{R_0}} \geq \tfrac{5}{8} r^* > \tfrac{1}{2} r^* 
\end{equation}
thus $\inf_{U\in \widehat{S}} \norm{u-\tilde{U}(\cdot-q)}_{B_{R_0}(q)}\geq \frac{1}{2} r^*$ and hence $d(q,u)=0$ for $q\notin B_{2R_0}(p)$; this means that $\supp(d(\cdot,u))\subset B_{2R_0}(p)$.

We show next that $d(\cdot, u)$ is equal to $1$ on a ball. Indeed if $u=U(\cdot-p)+\varphi$
\begin{eqnarray*}
\lefteqn{\inf_{\tilde{U}\in \widehat{S}} \norm{u-\tilde{U}(\cdot-q)}_{B_{R_0}(q)} \leq \norm{u-U(\cdot-q)}_{B_{R_0}(q)}}\\
& \leq& \norm{U(\cdot-p)-U(\cdot-q)}_{B_{R_0}(q)} + \norm{\varphi}_{B_{R_0}(q)} \leq \norm{\tau_{p-q}U-U}_{B_{R_0}} + \tfrac{1}{8}r^*.
\end{eqnarray*}
We can make the first term as small as we want by taking $|p-q|$ small, that is 
$$\inf_{U\in \widehat{S}} \norm{u-\tilde{U}(\cdot-q)}_{B_{R_0}(q)} \leq \tfrac{1}{4}r^*$$
 for $q\in B_r(p)$, $r$ small, which implies $d(q,u)=1$.

By the fact that $B_r(p)\subset \supp(d(\cdot, u))\subset B_{2R_0}(p)$ we have the well posedness of $\Upsilon(u)$ and
$$\Upsilon(u) = \frac{\int_{B_{2R_0}(p)} q \, d(q,u) dq}{\int_{B_{2R_0}(p)} d(q,u) dq}.$$
The main property comes straightforward, as well as the shift equivariance. We show now the continuity. 
Indeed, assume $\norm{u-v}_{H^s(\R^N)}\leq \frac{1}{8}r^*$. Then, by \eqref{eq_dim_58},
$$\norm{v-\tilde{U}(\cdot -q)}_{B_{R_0}(q)} \geq \norm{u-\tilde{U}(\cdot -q)}_{B_{R_0}(q)} - \norm{v-u}_{B_{R_0}(q)} \geq \tfrac{1}{2} r^*$$
and again we can conclude that $\supp(d(\cdot, v))\subset B_{2R_0}(p)$ for each $\norm{u-v}_{H^s(\R^N)}\leq \tfrac{1}{8} r^*$, where $p$ depends only on $u$. Moreover, observe that $\int_{B_{2R_0}(p)} d(q,u) dq \geq \int_{B_{r}(p)}1 \,dq \geq |B_{r}|=:C_1$ not depending on $u$ and $p$ (and similarly $C_2:= |B_{2R_0}|$), and that $d(q, \cdot)$ is Lipschitz 
(since $\psi$ and the norm are so, and the infimum over a family of Lipschitz function is still Lipschitz). 
Thus we have
\begin{eqnarray*}
\abs{\Upsilon(u)-\Upsilon(v)} &\leq& \frac{\int_{B_{2R_0}(p)} \abs{q} \, \abs{d(q,u)-d(q,v)} dq}{\int_{B_{2R_0}(p)} d(q,u) dq} +\\
&&+\int_{B_{2R_0}(p)} \abs{q} \, d(q,v) dq\frac{\int_{B_{2R_0}(p)} \abs{d(q,v)-d(q,u)} dq}{\int_{B_{2R_0}(p)} d(q,u) dq \int_{B_{2R_0}(p)} d(q,v) dq}\\
&\leq& \int_{B_{2R_0}(p)} \abs{q} dq \frac{1}{C_1}\left(1+ \frac{C_2}{C_1}\right)\norm{u-v}_{H^s(\R^N)}
=: C_p \norm{u-v}_{H^s(\R^N)}.
\end{eqnarray*}
Since $C_p$ can be bounded above by a constant of the type $C(1+\Upsilon(u))$, we have
$$\norm{u-v}_{H^s(\R^N)}\leq r_0 \implies \abs{\Upsilon(u)-\Upsilon(v)} \leq C(1+\Upsilon(u)) \norm{u-v}_{H^s(\R^N)};$$
in particular this implies the continuity.
\QED

\section{Singularly perturbed equation}
\label{sez_singur_perturb}

We come back now to our equation
\begin{equation}\label{eq_princ_epsx}
(-\Delta)^s u + V(\eps x) u= f(u), \quad x \in \R^N.
\end{equation}
It is known that the solutions of \eqref{eq_princ_epsx} can be characterized as critical points of the functional $I_{\eps}: H^s(\R^N) \to \R$
$$I_{\eps}(u):=\frac{1}{2} \norm{(-\Delta)^{s/2} u}_2^2 + \frac{1}{2} \int_{\R^N} V(\eps x) u^2 dx - \int_{\R^N}F(u)dx, \quad u \in H^s(\R^N) $$
where $I_{\eps}\in C^1(H^s(\R^N),\R)$, since $\norm{\cdot}_{H_{\eps}^s(\R^N)}$ is a norm.

We start with a technical result. Let $\nu_1$ be as in Lemma \ref{lem_tec2}; we want to show that the claim of the lemma 
continues holding, for $\eps$ small, if we replace $L_{m_0+\nu_1}$ with $I_{\eps}$, and $\widehat{S}$ with $S(r_0')\cap \{\eps \Upsilon(u) \in \Omega[\nu_1,\nu_0]\}$, $r_0'$ small.

\begin{lemma}\label{lemma_stimabasso_I}
Let $\nu_1$ and $\delta_0$ be as in Lemma \ref{lem_tec2}. 
Then there exist $\delta_1 \in (0, \delta_0)$ and $r_0'=r_0'(\delta_1) \in (0,r_0)$ sufficiently small, such that for every $\eps$ small we have 
$$I_{\eps}(u)\geq E_{m_0} + \delta_1$$
for each $u \in \{u \in S(r_0')\mid \eps \Upsilon(u) \in \Omega[\nu_1,+\infty)\} \supset \{u \in S(r_0')\mid \eps \Upsilon(u) \in \Omega[\nu_1,\nu_0]\}$.
\end{lemma}

\claim Proof.
First we improve Lemma \ref{lem_tec2} for $L_a$ in the direction of the nonautonomus equation. Indeed, by the assumption, we have $V(\eps \Upsilon(u))\geq m_0+\nu_1$, that is
$$L_{V(\eps \Upsilon(u))}(U) \geq L_{m_0+\nu_1}(U) \geq E_{m_0} + \delta_0$$
for any $U\in \widehat{S}$. Moreover, if $u=\tilde{U}(\cdot - p) + \tilde{\varphi}\in S(r_0)$ then, by Lemma \ref{lem_def_bar}, $\eps p \in \Omega_{2\eps R_0}\subset \Omega_{2R_0}$ which is compact. 
 By uniform continuity of $V$ and boundedness from above of $\widehat{S}$, we have
\begin{equation}\label{eq_dim_nonaut_to_aut}
L_{V(\eps p)}(U) \geq E_{m_0} + \delta_0/2
\end{equation}
for all $U\in \widehat{S}$ and $\eps$ small enough. 

Let now $r_0'$ to be fixed and $u=U(\cdot-p)+\varphi\in S(r_1)$. Then we have
$$I_{\eps}(u) = I_{\eps}(U(\cdot-p)+\varphi) = I_{\eps}(U(\cdot-p)) + I'_{\eps}(v)\varphi$$
for some $v\in H^s(\R^N)$ in the segment $[U(\cdot-p), u]$. Notice that $v$ lies in a ball of radius $\max \widehat{S}+r_0'$ and $I'_{\eps}$ sends bounded sets in bounded sets (uniformly on $\eps$); thus there exists a constant $C$, not depending on $U$, $p$ and $\varphi$, such that
\begin{equation}\label{eq_dim_stimaI}
I_{\eps}(u) \geq I_{\eps}(U(\cdot-p)) - C \norm{\varphi}_{H^s(\R^N)} \geq I_{\eps}(U(\cdot-p)) - \delta_1/2
\end{equation}
for $\norm{\varphi}_{H^s(\R^N)}< r_0'$ sufficiently small. 
Recalled that $\eps p \in \Omega_{2 R_0}$ we have, by the uniform continuity of $V$ and the uniform estimate \eqref{eq_stima_uniformeU}, for sufficiently small $\eps$, 
$$I_{\eps}(U(\cdot-p)) \geq L_{V(\eps p)}(U) - \delta_1/2$$
and the claim comes from \eqref{eq_dim_stimaI} and \eqref{eq_dim_nonaut_to_aut}, since, for $\delta_1< \delta_0/4$, 
$$I_{\eps}(u) \geq E_{m_0}+\delta_0/2- \delta_1 \geq E_{m_0} + \delta_1. $$
\QED

\bigskip

Before introducing the penalized functional, we state another technical lemma, which gives a (trivial, but useful) lower bound for $I'_{\eps}(v)v$ for small values of $v\in H^s(\R^N)$.

\begin{lemma}\label{lemma_stima_I}
There exists $r_1>0$ sufficiently small (see also the Remark below) and a constant $C>0$ such that
\begin{equation}\label{eq_stima_I}
I'_{\eps}(v)v \geq C\norm{v}_{H^s(\R^N)}^2
\end{equation}
for every $\eps>0$ and $v\in H^s(\R^N)$ with $\norm{v}_{H^s(\R^N)}\leq r_1$.
\end{lemma}

\claim Proof.
We have, by \eqref{eq_prop_f} with $\beta < \frac{1}{2}\underline{V}$,
\begin{eqnarray}
I'_{\eps}(v)v &\geq& \norm{(-\Delta)^{s/2}v}_2^2 + \int_{\R^N} \underline{V} v^2 dx - \beta \norm{v}_2^2 - C_{\beta} \norm{v}_{p+1}^{p+1} \notag \\
&\geq& \norm{(-\Delta)^{s/2}v}_2^2 + \frac{1}{2}\underline{V} \norm{v}_2^2 - C_{\beta} \norm{v}_{p+1}^{p+1} \label{eq_dim_da_migliorare}\\
&\geq& C\norm{v}_{H^s(\R^N)}^2 - C_{\beta} \norm{v}_{H^s(\R^N)}^{p+1}
\geq C'\norm{v}_{H^s(\R^N)}^2 \notag
\end{eqnarray}
where the last inequality holds for $\norm{v}_{H^s(\R^N)}$ small, since $p+1>2$.
\QED

\begin{remark} 
For a later use, we observe that one can improve \eqref{eq_dim_da_migliorare} by
\begin{equation}\label{eq_stima_2p}
\norm{(-\Delta)^{s/2}v}_2^2 + \frac{1}{2}\underline{V} \norm{v}_2^2 - 2^p C_{\beta} \norm{v}_{p+1}^{p+1} \geq C\norm{v}_{H^s(\R^N)}^2
\end{equation}
up to choosing a smaller $r_1$.
\end{remark}

\subsection{A penalized functional}

We want to study now a penalized functional (see \cite{CJT}), that is $I_{\eps}$ plus a term which forces solutions to stay in $\Omega$; we highlight that if we had information on $V$ at infinity (for example, trapping potentials) we would not have need of this penalization. 

Since $V>m_0$ on $\partial \Omega$, we can find an annulus around $\partial \Omega$ where this relation keeps holding, that is
$$V(x)>m_0, \quad \textnormal{ for $x\in \overline{\Omega_{2h_0}\setminus \Omega}$}$$
for $h_0$ sufficiently small. We then define the \emph{penalization} (or \emph{mass-concentrating}) functional $Q_{\eps}: H^s(\R^N) \to \R$
$$Q_{\eps}(u):=\left(\frac{1}{\eps^{\alpha}} \norm{u}^2_{L^2(\R^N\setminus (\Omega_{2h_0}/\eps))}-1\right)_+^{\frac{p+1}{2}}, \quad u \in H^s(\R^N)$$
where $\alpha \in (0, \min\{1/2, s\})$. 

We observe that, for every $u,v \in H^s(\R^N)$,
$$Q'_{\eps}(u)v = \frac{(p+1)}{\eps^{\alpha}} \left(\frac{1}{\eps^{\alpha}} \norm{u}^2_{L^2(\R^N\setminus (\Omega_{2h_0}/\eps))}-1\right)_+^{\frac{p-1}{2}} \int_{\R^N\setminus (\Omega_{2h_0}/\eps)} u v \, dx $$
and it is straightforward to prove the following estimate
\begin{equation} \label{eq_stima_Q}
Q'_{\eps}(u)u \geq (p+1) Q_{\eps}(u).
\end{equation}
We further set
$$J_{\eps}:=I_{\eps} + Q_{\eps}$$
the \emph{penalized} functional. It results that $Q_{\eps}$ and $J_{\eps}$ are in $C^1(H^s(\R^N), \R)$.

We want to find critical points of $J_{\eps}$ and show, afterwards, that these critical points, under suitable assumptions, are critical points of $I_{\eps}$ too, since $Q_{\eps}$ will be identically zero. 
Let $\eps=1$: observed that $Q_1(u)$ vanishes if $u$ have much mass inside $\Omega$, we see that $J_1(u)=I_1(u)$ holds when the mass of $u$ \emph{concentrates} in $\Omega$; this is why we say that $Q_{\eps}$ \emph{forces} $u$ to stay in $\Omega$. Similarly, as $\eps \to 0$, much less mass must be found outside $\Omega/\eps$.

We start by two technical lemmas. The first one gives a sufficient condition to pass from weak to strong convergent sequences in a Hilbert space, similarly to the convergence of the norms.

\begin{lemma}\label{lemma_conv_forte} 
Fix $\eps>0$ and let $(u_j)_j\subset H^s(\R^N)$ be such that
\begin{equation}\label{eq_ip_conv_J'}
\norm{J'_{\eps}(u_j)}_{H^{\text{-} s}(\R^N)} \to 0 \quad \textit{ as $j \to +\infty$}.
\end{equation}
Assume moreover that $u_j \wto u_0$ in $H^s(\R^N)$ as $j \to +\infty$, and that
\begin{equation}\label{eq_doppio_limite}
\lim_{R,j \to +\infty} \norm{u_j}_{L^q(\R^N \setminus B_R)}=0
\end{equation}
for $q=2$ and $q=p+1$. 
Then $u_j \to u_0$ in $H^s(\R^N)$ as $j \to +\infty$.
\end{lemma}

\claim Proof. 
We have by the weak lower semicontinuity of the norm
\begin{equation}\label{eq_semicont_inf}
\liminf_{j\to +\infty} \norm{u_j}_{H_{\eps}^s(\R^N)} \geq \norm{u_0}_{H_{\eps}^s(\R^N)}.
\end{equation}
Moreover
\begin{eqnarray*}
\lefteqn{\norm{u_j}_{H_{\eps}^s(\R^N)}^2 = \int_{\R^N} f(u_j) u_j dx + I'_{\eps}(u_j)u_j  } \\
&=& \Bigg( \int_{\R^N} f(u_j) u_j dx -\int_{\R^N} f(u_0)u_0 dx \Bigg) +\Big( I'_{\eps}(u_j)u_j -I'_{\eps}(u_0)u_0\Big) +\\
&&+ I'_{\eps}(u_0)u_0 + \int_{\R^N} f(u_0)u_0 dx =: (I) + (II) + \norm{u_0}_{H_{\eps}^s(\R^N)}^2;
\end{eqnarray*}
if we prove that
$$\limsup_{j \to +\infty} \big( (I) + (II) \big) \leq 0$$
we are done, because together with \eqref{eq_semicont_inf} we obtain
$$\norm{u_j}_{H_{\eps}^s(\R^N)} \to \norm{u_0}_{H_{\eps}^s(\R^N)} \quad \textnormal{ as $j \to +\infty$},$$
which implies the claim, since $H^s_{\eps}(\R^N)$ is a Hilbert space.

Focus on $(I)$; we have
\begin{eqnarray*}
\lefteqn{\int_{\R^N} (f(u_j)u_j - f(u_0)u_0) dx = \int_{B_R} (f(u_j)u_j - f(u_0)u_0) dx +}\\
&& + \int_{\R^N\setminus B_R} (f(u_j)u_j - f(u_0)u_0) dx =:(I_1) + (I_2).
\end{eqnarray*}
The piece $(I_2)$ can be made small for $j$ and $R$ sufficiently large, by exploiting the estimates on $f$, assumption \eqref{eq_doppio_limite} and the absolute continuity of the Lebesgue integral for $u_0$. 
For such large $R$ and $j$, up to taking a larger $j$, we can make the piece $(I_1)$ small by classical arguments. 

Focus now on $(II)$; we first observe that by exploiting H\"older inequalities and again classical arguments we have $I'_{\eps}(u_j)u_0 \to I'_{\eps}(u_0)u_0$.
Thus we have, by \eqref{eq_ip_conv_J'},
 \begin{eqnarray*}
\lefteqn{\limsup_{j \to +\infty} \left(I'_{\eps}(u_j) u_j - I'_{\eps}(u_0)u_0 \right)
= -\liminf_{j \to +\infty} \left(Q'_{\eps}(u_j) u_j - Q'_{\eps}(u_j)u_0 \right)} \\
&=& -\liminf_{j \to +\infty} \Bigg( \left(\frac{1}{\eps^{\alpha}} \norm{u_j}^2_{L^2(\R^N\setminus (\Omega_{2h_0}/\eps))}-1\right)_+^{\frac{p-1}{2}} \cdot \\
&& \cdot \Bigg(\int_{\R^N\setminus (\Omega_{2h_0}/\eps)}u_j^2 dx -\int_{\R^N\setminus (\Omega_{2h_0}/\eps)} u_j u_0 dx \Bigg) \Bigg) \leq 0
 \end{eqnarray*}
where the last inequality is due to the following fact: observe first that $u_j \wto u_0$ in $H^s_{\eps}(\R^N)\hookrightarrow L^2(\R^N)$ thus (by restriction) $u_j \wto u_0$ in $L^2(\R^n \setminus (\Omega_{2h_0}/\eps))$; by definition of weak convergence and by the lower semicontinuity of the norm, we have
$$\liminf_{j \to +\infty} \int_{\R^N\setminus (\Omega_{2h_0}/\eps)}u_j^2 dx \geq \int_{\R^N\setminus (\Omega_{2h_0}/\eps)}u_0^2 dx = \lim_{j \to +\infty} \int_{\R^N\setminus (\Omega_{2h_0}/\eps)}u_j u_0 dx ,$$
that is
$$\liminf_{j \to +\infty} \Bigg(\int_{\R^N\setminus (\Omega_{2h_0}/\eps)}u_j^2dx -\int_{\R^N\setminus (\Omega_{2h_0}/\eps)} u_j u_0 dx \Bigg) \geq 0.$$
Noticed that $a_n \geq 0$ and $\liminf_n b_n \geq 0$ imply $\liminf_n (a_n b_n) \geq 0$, we conclude.
\QED

\bigskip

The second Lemma is a lower bound for $J_{\eps}$ with respect to the functional $L_{m_0}$.
We highlight that in what follows we understand that the case $m_0=\underline{V}$, i.e. $m_0$ \emph{global} minimum, gives rise to a not-perturbed result.

\begin{lemma}\label{lemma_stima_J}
Set $C_{min}:= \frac{1}{2}(m_0-\underline{V})\geq 0$ we have, for $\eps$ small and $u\in H^s(\R^N)$,
$$J_{\eps}(u) \geq L_{m_0}(u) - C_{min} \eps^{\alpha}.$$
\end{lemma}

\claim Proof.
We have, recalling that $m_0$ is the infimum of $V$ over $\Omega_{2h_0}$ and $\underline{V}$ is the infimum over $\R^N$,
\begin{eqnarray*}
J_{\eps}(u) &=& L_{m_0}(u) + \frac{1}{2} \int_{\R^N} (V(\eps x) - m_0) u^2 dx + Q_{\eps}(u)\\
&\geq & L_{m_0}(u) + \frac{1}{2} \int_{\R^N\setminus(\Omega_{2h_0} / \eps)} (V(\eps x) - m_0) u^2 dx + Q_{\eps}(u)\\
&\geq &
 L_{m_0}(u) - C_{min} \norm{u}_{L^2(\R^N\setminus (\Omega_{2h_0}/\eps))}^2 + Q_{\eps}(u).
\end{eqnarray*}
If $\norm{u}_{L^2(\R^N\setminus (\Omega_{2h_0}/\eps))}^2 \leq 2\eps^{\alpha}$ we have the claim by the positivity of $Q_{\eps}(u)$. If instead $\norm{u}_{L^2(\R^N\setminus (\Omega_{2h_0}/\eps))}^2 \geq 2\eps^{\alpha}$, then
$$ Q_{\eps}(u) \geq \left(\frac{1}{2\eps^{\alpha}} \norm{u}^2_{L^2(\R^N\setminus (\Omega_{2h_0}/\eps))}\right)^{\frac{p+1}{2}}
\geq \frac{1}{2\eps^{\alpha}} \norm{u}^2_{L^2(\R^N\setminus (\Omega_{2h_0}/\eps))}. $$
Thus
$$ J_{\eps}(u) \geq L_{m_0} + \left (\frac{1}{2\eps^{\alpha}} - C_{min} \right) \norm{u}_{L^2(\R^N\setminus (\Omega_{2h_0}/\eps))}^2 
\geq L_{m_0} $$
for $\eps$ small. This concludes the proof.
\QED

\subsection{Critical points and Palais-Smale condition}

In order to get critical points of $J_{\eps}$ we want to implement a deformation argument. As usual, we need a uniform estimate from below of $\norm{J'_{\eps}(u)}_{H^{\text{-}s}(\R^N)}$, and this is the next goal. 

First, by the strict monotonicity of $E_a$, let us fix $l_0'=l_0'(\nu_1)>0$ such that
$$E_{m_0}< l_0' < E_{m_0 + \nu_1};$$
as well as $\nu_0$ and $l_0$, even $l_0'$ will be let vary as $(l_0')^n \to E_{m_0}$ in the proof of the existence.

\begin{lemma}\label{lemma_stima_basso} 
Let $r_0$ and $r_1$ be as in Lemma \ref{lem_def_bar} and Lemma \ref{lemma_stima_I}. There exists $r_2'\in (0,\min\{r_0, r_1\})$ sufficiently small with the following property: let $0<\rho_1<\rho_0\leq r_2'$ and $(u_{\eps})_{\eps}\subset S(r_2')$ be such that
\begin{eqnarray}
&\quad \norm{J'_{\eps}(u_{\eps})}_{H^{\text{-}s}(\R^N)}\to 0 \quad \textit{ as $\eps \to 0$},& \label{eq_PS2}\\ 
&J_{\eps}(u_{\eps})\leq l_0' < E_{m_0+\nu_1}, \quad \textit{ for any $\eps >0$},& \label{eq_PS1} 
\end{eqnarray}
with the additional assumption
$$(\widehat{\rho}(u_{\eps}))_{\eps} \subset [0,\rho_0], \quad (\eps\Upsilon(u_{\eps}))_{\eps} \subset \Omega[0,\nu_0].$$
Then, for $\eps$ small
$$\widehat{\rho}(u_{\eps})\in [0,\rho_1], \quad \eps\Upsilon(u_{\eps})\in \Omega[0,\nu_1].$$
\end{lemma}

We notice, by \eqref{eq_PS2} and \eqref{eq_PS1}, that $(u_{\eps})_{\eps}$ resembles a particular (truncated) Palais-Smale sequence. As an immediate consequence of the Lemma, set the sublevel
$$J_{\eps}^c:= \{ u \in H^s(\R^N) \mid J_{\eps}(u) \leq c\}$$
 we have the following theorem.

\begin{theorem}\label{theorem_stima_J'}
There exists $r_2'\in (0,\min\{r_0, r_1\})$ sufficiently small with the following property: if $0<\rho_1<\rho_0\leq r_2'$, then there exists a $\delta_2=\delta_2(\rho_0,\rho_1)>0$ such that, for $\eps$ small
$$\norm{J'_{\eps}(u)}_{H^{\text{-}s}(\R^N)}\geq \delta_2$$
for any 
\begin{eqnarray*}
u &\in& \left\{u\in S(r_2')\cap J_{\eps}^{l_0'} \mid (\widehat{\rho}(u), \eps \Upsilon(u))\in \big((0,\rho_0]\times \Omega[0,\nu_0]\big)\setminus\big([0,\rho_1]\times \Omega[0,\nu_1]\big)\right\} \\
&\supset& \left\{u\in S(r_2')\cap J_{\eps}^{l_0'} \mid \rho_1<\widehat{\rho}(u)\leq \rho_0, \; \eps \Upsilon(u)\in \Omega(\nu_1,\nu_0]\right\}.
\end{eqnarray*}
\end{theorem}

\begin{remark}\label{remark_triplanorma}
Arguing as in the last part of the proof of Lemma \ref{lemma_stima_unif}, noticed that $\widehat{S}$ is compact not only in $H^s(\R^N)$ but also in $L^q(\R^N)$ for $q \in [2, 2^*_s]$, if $(U_n)_n\subset \widehat{S}$ and $(\theta_n)_n \subset \R^N$ is included in a compact, we have
$$\lim_{n\to +\infty} \tnorm{U_n(\cdot + \theta_n)}_{\R^N \setminus B_n}=0,$$
where the norm $\tnorm{\cdot}_{\R^N\setminus B_n}$ is defined in \eqref{eq_def_triplanorm}.
\end{remark}

\claim Proof of Lemma \ref{lemma_stima_basso}. 
We use the notation, for $h>0$,
$$\Omega^{\eps}_{h}:=(\Omega_{\eps h})/\eps = (\Omega/\eps)_h$$
and notice that if $h<h'$ then $\Omega/\eps \subset\Omega^{\eps}_{h} \subset \Omega^{\eps}_{h'}$. 
Let $r_2'<\min\{r_0, r_1\}$ to be fixed.

\medskip

\noindent \textbf{Step 1.} \textit{An estimate for $u_{\eps}$.}

We have, for $u_{\eps}= U_{\eps}(\cdot-p_{\eps}) + \varphi_{\eps} \in S(r_2')$, 
\begin{eqnarray*}
\tnorm{u_{\eps}}_{\R^N\setminus (\Omega/\eps)} &\leq& 
\tnorm{U_{\eps}(\cdot-p_{\eps})}_{\R^N\setminus (\Omega/\eps)} + \tnorm{\varphi_{\eps}}_{\R^N\setminus (\Omega/\eps)}\\
&\leq& \tnorm{U_{\eps}(\cdot-p_{\eps}+\Upsilon(u_{\eps}))}_{\R^N\setminus (\Omega/\eps-\Upsilon(u_{\eps}))} + C
\norm{\varphi_{\eps}}_{H^s(\R^N)}\\
&\leq& \tnorm{U_{\eps}(\cdot-p_{\eps}+\Upsilon(u_{\eps}))}_{\R^N\setminus (\Omega/\eps-\Upsilon(u_{\eps}))} + C
 r_2'.
\end{eqnarray*}
By the fact that $\eps \Upsilon(u_{\eps})\in \Omega[\nu_1,\nu_0]\subset \Omega$, we have that $0 \in \Omega/\eps - \Upsilon(u_{\eps})$ and thus $\Omega/\eps - \Upsilon(u_{\eps})$ expands in $\R^N$ as $\eps \to 0$. Moreover by Lemma \ref{lem_def_bar} we have $\theta_{\eps}:= \Upsilon(u_{\eps})- p_{\eps} \in B_{2R_0}$ compact. By Remark \ref{remark_triplanorma}, for $\varepsilon$ small we have
\begin{equation} \label{eq_stima_princ_ueps}
\tnorm{u_{\eps}}_{\R^N\setminus (\Omega/\eps)}\leq (1+C)r_2' =C' r_2'
\end{equation}
Let 
$$n_{\eps}:= \left[ \frac{\sqrt{1+4h_0/\eps}+1}{2}\right]\in \N$$
which by definition satisfies $\eps n_{\eps}(n_{\eps}+1)\leq h_0$ and $n_{\eps}\to+ \infty$ as $\eps \to 0$. We have
$$\sum_{i=1}^{n_{\eps}} \norm{u_{\eps}}_{L^2(\Omega_{n_{\eps}(i+1)}^{\eps}\setminus \Omega_{n_{\eps} i}^{\eps})}^2 
\leq 
 \norm{u_{\eps}}_{L^2(\Omega_{n_{\eps}(n_{\eps}+1)}^{\eps}\setminus \Omega_{n_{\eps}}^{\eps})}^2 \leq (C' r_2')^2$$
and similarly 
$$\sum_{i=1}^{n_{\eps}} [u_{\eps}]_{\Omega_{n_{\eps}(i+1)}^{\eps}\setminus \Omega_{n_{\eps} i}^{\eps}, \R^N}^2 \leq (C'
 r_2')^2, \qquad \sum_{i=1}^{n_{\eps}} \norm{u_{\eps}}_{L^{p+1}(\Omega_{n_{\eps}(i+1)}^{\eps}\setminus \Omega_{n_{\eps} i}^{\eps})}^{p+1} \leq (C'
 r_2')^{p+1}$$
thus, for some $C=C(r_2')$,
$$\sum_{i=1}^{n_{\eps}} \left( 
\norm{u_{\eps}}_{\Omega_{n_{\eps}(i+1)}^{\eps}\setminus \Omega_{n_{\eps} i}^{\eps}}^2 + \norm{u_{\eps}}_{L^{p+1}(\Omega_{n_{\eps}(i+1)}^{\eps}\setminus \Omega_{n_{\eps} i}^{\eps})}^{p+1} \right) \leq C. $$
This implies that there exists $i_{\eps} \in \{1, \dots, n_{\eps}\}$ such that
\begin{equation}\label{eq_stima_anello}
\norm{u_{\eps}}_{A^{\eps}}^2 + \norm{u_{\eps}}_{L^{p+1}(A^{\eps})}^{p+1} \leq \frac{C}{n_{\eps}}\to 0 \quad \textnormal{ as $\eps \to 0$},
\end{equation}
where 
$$A^{\eps}:= \Omega_{n_{\eps}(i_{\eps}+1)}^{\eps}\setminus \Omega_{n_{\eps} i_{\eps}}^{\eps}$$
and $C$ depends on $r_2'$ (we will omit this dependence).

\medskip

\noindent \textbf{Step 2.} \textit{Split the sequence.}

Consider cutoff functions $\varphi_{\eps} \in C^{\infty}_c(\R^N)$ 
$$\Omega^{\eps}_{n_{\eps} i_{\eps}} \prec \varphi_{\eps} \prec \Omega^{\eps}_{n_{\eps}(i_{\eps}+1)}$$
such that $\norm{\nabla \varphi_{\eps}}_{\infty} \leq \frac{C}{n_{\eps}}= o(1)$ as $\eps \to 0$ (which is possible because the distance between $\Omega^{\eps}_{n_{\eps} i_{\eps}}$ and $\Omega^{\eps}_{n_{\eps}(i_{\eps}+1)}$ is $n_{\eps}\to +\infty$).

Define 
$$u_{\eps}^{(1)}:=\varphi_{\eps} u_{\eps}, \quad u_{\eps}^{(2)}:=(1-\varphi_{\eps}) u_{\eps}\quad \textnormal{ and } \; u_{\eps}=u_{\eps}^{(1)}+u_{\eps}^{(2)};$$
notice that both $\supp \big(u_{\eps}^{(1)} u_{\eps}^{(2)}\big)$ and $\supp(F(u_{\eps}) -F(u_{\eps}^{(1)}) - F(u_{\eps}^{(2)}))$ are contained in $A^{\eps}$, that is where we gained the estimate of the norm. Moreover, since $$
\supp(u_{\eps}^{(1)}) \subset \Omega_{n_{\eps}(i_{\eps}+1)}^{\eps} \subset \Omega_{\eps n_{\eps}(n_{\eps}+1)}/\eps \subset \Omega_{2h_0}/\eps $$
we have, by definition of $Q_{\eps}$, that $Q_{\eps}(u_{\eps}^{(1)})=0$, $Q_{\eps}(u_{\eps}) = Q_{\eps}(u_{\eps}^{(2)})$ and
\begin{equation}\label{eq_relaz_der_Qeps}
Q'_{\eps}(u_{\eps}^{(1)})=0, \quad Q'_{\eps}(u_{\eps}) = Q'_{\eps}(u_{\eps}^{(2)}).
\end{equation}

\medskip

\noindent \textbf{Step 3.} \textit{Relations of the functionals.}

We show that 
$$\abs{I_{\eps}(u_{\eps}) - I_{\eps}(u_{\eps}^{(1)}) - I_{\eps}(u_{\eps}^{(2)})}\to 0 \quad \textnormal{ as $\eps \to 0$}$$ 
from which
\begin{equation}\label{eq_relaz_funzionali}
J_{\eps}(u_{\eps}) = I_{\eps}(u_{\eps}^{(1)}) + I_{\eps}(u_{\eps}^{(2)}) + Q_{\eps}(u_{\eps}^{(2)}) + o(1).
\end{equation}
Indeed 
\begin{eqnarray*}
\lefteqn{\abs{I_{\eps}(u_{\eps}) - I_{\eps}(u_{\eps}^{(1)}) - I_{\eps}(u_{\eps}^{(2)})}
 \leq \Bigg|\int_{ \R^N} (-\Delta)^{s/2}u_{\eps}^{(1)} (-\Delta)^{s/2}u_{\eps}^{(2)}dx \Bigg| +}\\
 &&+ \int_{ A^{\eps}} \abs{V(\eps x)u_{\eps}^{(1)} u_{\eps}^{(2)}}dx +\int_{ A^{\eps}} \abs{F(u_{\eps}) -F(u_{\eps}^{(1)}) - F(u_{\eps}^{(2)}) }dx \\
&=:& (I)+(II)+(III).
\end{eqnarray*}

The second piece can be easily estimated by the boundedness of $\varphi_{\eps}$ and $V$, and the information on the $2$-norm given by \eqref{eq_stima_anello}, i.e $(II) \leq \frac{C}{n_{\eps}}$. 
Similarly, as regards $(III)$, we estimate each single piece separately, in the same way: use \eqref{eq_prop_f} and the information on the $2$-norm and the $(p+1)$-norm given by \eqref{eq_stima_anello}, obtaining $(III)\leq \frac{C}{n_{\eps}}$. 

Focus instead on (I). Recall that $(u_{\eps})_{\eps}\subset S(r_2')$, and thus $\norm{u_{\eps}}_2$ is bounded. We have 
\begin{eqnarray*}
(I) &\leq& C \int_{\R^{2N}} \frac{\abs{u_{\eps}^{(1)}(x)-u_{\eps}^{(1)}(y)} \abs{u_{\eps}^{(2)}(x)-u_{\eps}^{(2)}(y)}}{|x-y|^{N+2s}} dx \, dy \\
&\leq& 2C \int_{\Omega^{\eps}_{n_{\eps} i_{\eps}}\times \complement(\Omega^{\eps}_{n_{\eps}(i_{\eps}+1)})} \frac{\abs{u_{\eps}^{(1)}(x)-u_{\eps}^{(1)}(y)} \abs{u_{\eps}^{(2)}(x)-u_{\eps}^{(2)}(y)}}{|x-y|^{N+2s}} dx \, dy+ \\
&&+ 2C \int_{A^{\eps}\times \R^N} \frac{\abs{u_{\eps}^{(1)}(x)-u_{\eps}^{(1)}(y)} \abs{u_{\eps}^{(2)}(x)-u_{\eps}^{(2)}(y)}}{|x-y|^{N+2s}} dx \, dy \\
&=:& 2C \big((I_1) + (I_2)\big)
\end{eqnarray*}
since on $\Omega^{\eps}_{n_{\eps} i_{\eps}} \times \Omega^{\eps}_{n_{\eps} i_{\eps}}$ and $\complement(\Omega^{\eps}_{n_{\eps}(i_{\eps}+1)}) \times \complement(\Omega^{\eps}_{n_{\eps}(i_{\eps}+1)})$ the integrand is null. Focusing on $(I_1)$
\begin{eqnarray*}
(I_1) &=&\int_{(\Omega^{\eps}_{n_{\eps} i_{\eps}}\times \complement(\Omega^{\eps}_{n_{\eps}(i_{\eps}+1)}))\cap \{|x-y|>n_{\eps}\}} \frac{\abs{u_{\eps}(x)u_{\eps}(y)}}{|x-y|^{N+2s}} dx \, dy\\
&\leq& \frac{1}{2}\int_{(\Omega^{\eps}_{n_{\eps} i_{\eps}}\times \complement(\Omega^{\eps}_{n_{\eps}(i_{\eps}+1)}))\cap \{|x-y|>n_{\eps}\}} \frac{u^2_{\eps}(x) + u^2_{\eps}(y)}{|x-y|^{N+2s}}dx \, dy \\
&=& \frac{1}{2} \int_{\Omega^{\eps}_{n_{\eps} i_{\eps}}}u^2_{\eps}(x) \int_{\complement(\Omega^{\eps}_{n_{\eps}(i_{\eps}+1)})\cap \{|x-y|>n_{\eps}\}} \frac{1}{|x-y|^{N+2s}} dy \, dx+\\
&&+ \frac{1}{2} \int_{\complement(\Omega^{\eps}_{n_{\eps}(i_{\eps}+1)})}u^2_{\eps}(y) \int_{\Omega^{\eps}_{n_{\eps} i_{\eps}}\cap \{|x-y|>n_{\eps}\}} \frac{1}{|x-y|^{N+2s}} dx \, dy \\
&\leq& C\norm{u_{\eps}}_2^2 \int_{|x-y|>n_{\eps}} \frac{1}{|x-y|^{N+2s}} dx \, dy 
\leq \frac{C}{{n_{\eps}}^{2s}} \to 0 \quad \textnormal{ as $\eps \to 0$}.
\end{eqnarray*}
Focusing on $(I_2)$ we have
\begin{eqnarray*}
\lefteqn{(I_2) }\\
&\leq& \int_{A^{\eps}\times \R^N} \frac{1}{|x-y|^{N+2s}} \Big( \abs{(\varphi_{\eps}(x) -\varphi_{\eps}(y))u_{\eps}(x)} + \abs{ \varphi_{\eps}(y)(u_{\eps}(x)-u_{\eps}(y))} \Big) \cdot \\
&& \cdot \Big( \abs{(\varphi_{\eps}(y)-\varphi_{\eps}(x))u_{\eps}(x)} + \abs{(1-\varphi_{\eps}(y))(u_{\eps}(x) -u_{\eps}(y))} \Big) dx \, dy \\
&\leq& \int_{A^{\eps}\times \R^N} \frac{\abs{\varphi_{\eps}(x) -\varphi_{\eps}(y)}^2 \abs{u_{\eps}(x)}^2}{|x-y|^{N+2s}} dx \, dy
 + \int_{A^{\eps}\times \R^N} \frac{\abs{u_{\eps}(x) -u_{\eps}(y)}^2}{|x-y|^{N+2s}}dx \, dy+\\
&&+ 2\int_{A^{\eps}\times \R^N} \frac{\abs{\varphi_{\eps}(x) -\varphi_{\eps}(y)} \abs{u_{\eps}(x)} \abs{u_{\eps}(x) -u_{\eps}(y)}}{|x-y|^{N+2s}}dx \, dy\\
&\leq& \int_{A^{\eps}\times \R^N} \frac{\abs{\varphi_{\eps}(x) -\varphi_{\eps}(y)}^2 \abs{u_{\eps}(x)}^2}{|x-y|^{N+2s}} dx \, dy+\int_{A^{\eps}\times \R^N} \frac{\abs{u_{\eps}(x) -u_{\eps}(y)}^2}{|x-y|^{N+2s}}dx \, dy+\\
&& +2 \Bigg( \int_{A^{\eps}\times \R^N} \frac{\abs{\varphi_{\eps}(x) -\varphi_{\eps}(y)}^2 \abs{u_{\eps}(x)}^2}{|x-y|^{N+2s}} dx \, dy \Bigg)^{\frac{1}{2}} \Bigg( \int_{A^{\eps}\times \R^N} \frac{\abs{u_{\eps}(x) -u_{\eps}(y)}^2}{|x-y|^{N+2s}}dx \, dy \Bigg)^{\frac{1}{2}} \\
&=:& A^2 + B^2 + 2AB
\end{eqnarray*}
and we see that both $A$ and $B$ go to zero: $B=[u_{\eps}]_{A^{\eps}, \R^N} \leq \frac{C}{n_{\eps}^{1/2}}$ by \eqref{eq_stima_anello}, while for $A$ we exploit that $\norm{\nabla \varphi_{\eps}}_{\infty}\to 0$. 
Indeed, let $\alpha_{\eps}:=\frac{1}{\norm{\nabla \varphi_{\eps}}_{\infty}}$; we have

\begin{eqnarray*}
A^2 &=& \int_{A^{\eps}} \abs{u_{\eps}(x)}^2 \int_{|x-y|\leq \alpha_{\eps}} \frac{\abs{\varphi_{\eps}(x) -\varphi_{\eps}(y)}^2}{|x-y|^{N+2s}} dy \, dx+\\
 &&+ \int_{A^{\eps}} \abs{u_{\eps}(x)}^2 \int_{|x-y|> \alpha_{\eps}} \frac{\abs{\varphi_{\eps}(x) -\varphi_{\eps}(y)}^2}{|x-y|^{N+2s}} dy \, dx\\
&\leq & C \norm{u_{\eps}}_{L^2(A^{\eps})}^2 \Bigg(\norm{\nabla \varphi_{\eps}}_{\infty}^2 \int_{|z|\leq \alpha_{\eps}} \frac{1}{|z|^{N+2s-2}} dz + 4 \norm{\varphi_{\eps}}_{\infty}^2 \int_{|z|> \alpha_{\eps}} \frac{1}{|z|^{N+2s}}dz \Bigg)\\
&\leq& \frac{C}{n_{\eps}} \alpha_{\eps}^{-2s} \left(\norm{\nabla \varphi_{\eps}}_{\infty}^2 \alpha_{\eps}^2 + 1\right)
= \frac{C}{n_{\eps}} \norm{\nabla \varphi_{\eps}}_{\infty}^{2s}\leq \frac{C}{n_{\eps}^{2s+1}}\to 0 \quad \textnormal{ as $\eps \to 0$}.
\end{eqnarray*}
Thus $(I_2) \leq \frac{C}{n_{\eps}^{2s+1}} + \frac{C}{n_{\eps}^{s+1}} + \frac{C}{n_{\eps}}\leq \frac{C'}{n_{\eps}}\to 0$, 
which reaches the claim.

\medskip

\noindent \textbf{Step 4.} \textit{Relations of the derivatives.}

We have
\begin{equation}\label{eq_conv_I'}
\norm{I'_{\eps}(u_{\eps}) - I'_{\eps}(u_{\eps}^{(1)}) - I'_{\eps}(u_{\eps}^{(2)})}_{H^{\text{-}s}(\R^N)}\to 0 \quad \textnormal{ as $\eps \to 0$},
\end{equation}
from which, joined to \eqref{eq_relaz_der_Qeps},
\begin{equation}\label{eq_relaz_derivate}
J'_{\eps}(u_{\eps} ) = I'_{\eps}(u_{\eps}^{(1)}) + I'_{\eps}(u_{\eps}^{(2)}) + Q'_{\eps}(u_{\eps}^{(2)}) + o(1).
\end{equation}
Indeed by H\"older inequality, for any $v\in H^s(\R^N)$,
$$ \abs{I'_{\eps}(u_{\eps})v - I'_{\eps}(u_{\eps}^{(1)})v - I'_{\eps}(u_{\eps}^{(2)})v} \leq \int_{A^{\eps}} \abs{f(u_{\eps})-f(u_{\eps}^{(1)}) - f(u_{\eps}^{(2)})} |v| dx $$
and again we argue in the same way as in the third piece of Step 3, observing that, by \eqref{eq_prop_f}, $|f(u_{\eps})| |v| \leq \beta |u_{\eps}| |v| + C_{\beta} |u_{\eps}|^{p} |v|$ thus
\begin{eqnarray*} 
\int_{A^{\eps}} |f(u_{\eps})| |v| dx &\leq& \beta \norm{u_{\eps}}_{L^2(A^{\eps})} \norm{v}_2 + C_{\beta}\norm{u_{\eps}}_{L^{p+1}(A^{\eps})}^p \norm{v}_{p+1} \\
&\leq& C\left( \beta \norm{u_{\eps}}_{L^2(A^{\eps})} + C_{\beta}\norm{u_{\eps}}_{L^{p+1}(A^{\eps})}^p \right) \norm{v}_{H^s}
\end{eqnarray*}
and hence the claim. In particular, $\abs{\big(I'_{\eps}(u_{\eps}) - I'_{\eps}(u_{\eps}^{(1)}) - I'_{\eps}(u_{\eps}^{(2)})\big)u_{\eps}^{(2)}} \leq \frac{C}{n_{\eps}}$.
We see also that 
\begin{equation}\label{eq_o_picc}
 I'_{\eps}(u_{\eps}^{(1)})u_{\eps}^{(2)} = o(1).
\end{equation}
Indeed
\begin{eqnarray*}
\lefteqn{\abs{I'_{\eps}(u_{\eps}^{(1)})u_{\eps}^{(2)}} \leq \Bigg|\int_{ \R^N} (-\Delta)^{s/2}u_{\eps}^{(1)} (-\Delta)^{s/2}u_{\eps}^{(2)} dx\Bigg|+ } \\
&&+ \int_{ A^{\eps}} \abs{V(\eps x)u_{\eps}^{(1)} u_{\eps}^{(2)}dx } + \int_{A^{\eps}} \abs{f(u_{\eps}^{(1)})u_{\eps}^{(2)}dx }
=: (I) + (II) + (III)
\end{eqnarray*}
where for $(I)$ and $(II)$ we argue as in Step 3 obtaining $(I)+ (II) \leq \frac{C}{n_{\eps}^{2s}} + \frac{C}{n_{\eps}}$, while for $(III)$ we argue as in \eqref{eq_relaz_derivate} obtaining $(III) \leq \frac{C}{n_{\eps}} $.

\medskip

\noindent \textbf{Step 5.} \textit{Convergence of $u_{\eps}^{(2)}$.}

Observing that the support of $u_{\eps}^{(2)}$ is outside $\Omega/\eps$, we have with arguments similar to Step 3 that, by \eqref{eq_stima_princ_ueps},
\begin{equation}\label{eq_dim_stima_r2}
\norm{u_{\eps}^{(2)}}_{H^s(\R^N)}\leq r_1.
\end{equation}
Indeed, focusing only on the nonlocal part, we have (recall that $\supp(u_{\eps}^{(2)})\subset \complement(\Omega/\eps)$)
\begin{eqnarray*}
\int_{\R^{2N}} \frac{\abs{u_{\eps}^{(2)}(x)-u_{\eps}^{(2)}(y)}^2}{|x-y|^{N+2s}}
&\leq& 2 \int_{\complement(\Omega/\eps)\times \R^N} \frac{\abs{u_{\eps}^{(2)}(x)-u_{\eps}^{(2)}(y)}^2}{|x-y|^{N+2s}} dx \, dy \\
&\leq & 4 \int_{\complement(\Omega/\eps)\times \R^{N}} \frac{\abs{\varphi_{\eps}(x)-\varphi_{\eps}(y)}^2\abs{u_{\eps}(x)}^2}{|x-y|^{N+2s}}dx \, dy+\\
&&+ 4 \int_{\complement(\Omega/\eps)\times \R^N} \frac{\abs{u_{\eps}(x)-u_{\eps}(y)}^2}{|x-y|^{N+2s}}dx \, dy
\end{eqnarray*}
and we use again the final argument in Step 3 and \eqref{eq_stima_princ_ueps} to gain, for $\eps $ small,
$$\norm{u_{\eps}^{(2)}}_{H^s(\R^N)}^2 \leq(C+ o(1)) \norm{u_{\eps}}_{L^2(\complement(\Omega/\eps))}^2 + C [u_{\eps}]_{\complement(\Omega/\eps), \R^N}^2 \leq (Cr_2')^2$$
where $C$ does not depend on $r_2'$. We choose thus $r_2'$ such that \eqref{eq_dim_stima_r2} holds.

This allows us to use Lemma \ref{lemma_stima_I}. 
By joining \eqref{eq_PS2}, \eqref{eq_relaz_derivate}, \eqref{eq_o_picc}, \eqref{eq_stima_I}, \eqref{eq_stima_Q} we obtain 
\begin{eqnarray*}
&o(1) = J'_{\eps}(u_{\eps} )u_{\eps}^{(2)} = I'_{\eps}(u_{\eps}^{(1)}) u_{\eps}^{(2)} + I'_{\eps}(u_{\eps}^{(2)})u_{\eps}^{(2)} + Q'_{\eps}(u_{\eps}^{(2)})u_{\eps}^{(2)} + o(1)& \\
&\geq C \norm{u_{\eps}^{(2)}}_{H^s(\R^N)}^2 + (p+1)Q_{\eps}(u_{\eps}^{(2)}) +o(1)&
\end{eqnarray*}
or more precisely (we highlight this for a later use), for some $C=C(r_2')$,
\begin{equation}
 o(1) = J'_{\eps}(u_{\eps} )u_{\eps}^{(2)} \geq C \norm{u_{\eps}^{(2)}}_{H^s(\R^N)}^2 + (p+1)Q_{\eps}(u_{\eps}^{(2)}) -
 \left( \frac{C}{n_{\eps}}+\frac{C}{n_{\eps}^{2s}}\right), \label{eq_dim_stima_mista}
\end{equation}
which implies (since $Q_{\eps}$ is positive) that
\begin{equation}\label{eq_conv_u2}
\norm{u_{\eps}^{(2)}}_{H^s(\R^N)}\to 0 \quad \textnormal{ as $\eps\to 0$}
\end{equation}
and
\begin{equation}\label{eq_conv_Q2}
Q_{\eps}(u_{\eps}^{(2)})\to 0 \quad \textnormal{ as $\eps\to 0$}.
\end{equation}
As a further consequence, \eqref{eq_conv_u2} and the boundedness of $V$ imply 
\begin{equation}\label{eq_conv_I2}
I_{\eps}(u_{\eps}^{(2)})\to 0 \quad \textnormal{and} \quad I'_{\eps}(u_{\eps}^{(2)})\to 0 \quad \textnormal{ as $\eps \to 0$}.
\end{equation}

\medskip

\noindent \textbf{Step 6.} \textit{Convergence of $I_{\eps}'(u_{\eps}^{(1)})$.}

In particular we obtain from \eqref{eq_conv_I2}, together with \eqref{eq_relaz_funzionali} and \eqref{eq_conv_Q2}, that
\begin{equation}\label{eq_confr_JI}
J_{\eps}(u_{\eps}) = I_{\eps}(u_{\eps}^{(1)}) + o(1).
\end{equation}
We want now to show that
\begin{equation}\label{eq_conv_I1}
I'_{\eps}(u_{\eps}^{(1)})\to 0 \quad \textnormal{in $H^{\text{-}s}(\R^N)$} \quad \textnormal{ as $\eps \to 0$}.
\end{equation}
Start observing that \eqref{eq_conv_I2} together with \eqref{eq_conv_I'} give
\begin{equation}\label{eq_info_I1}
I'_{\eps}(u_{\eps}) = I'_{\eps}(u_{\eps}^{(1)}) + o(1);
\end{equation}
let now $v \in H^s(\R^N)$ and evaluate $I'_{\eps}(u_{\eps}^{(1)})v$. We want to exploit \eqref{eq_info_I1} together again with the assumption \eqref{eq_PS2}. In order to do this we need to pass from $u_{\eps}^{(2)}$ to $u_{\eps}$, but getting rid of $Q'_{\eps}(u_{\eps})$ on which we have no information. Thus we introduce a cutoff function $\tilde{\varphi}\in C^{\infty}_c(\R^N)$ such that
$$\Omega_{\frac{3}{2}h_0} \prec \tilde{\varphi} \prec \Omega_{2h_0}$$
and hence
$$\supp(u_{\eps}^{(1)}) \subset \Omega_{h_0}/\eps\subset \Omega_{\frac{3}{2}h_0}/\eps\subset \{\tilde{\varphi}(\eps \cdot)\equiv 1\}\subset \supp(\tilde{\varphi}(\eps \cdot)) \subset \Omega_{2h_0}/\eps.$$
Thus we have
\begin{eqnarray*}
\lefteqn{ I'_{\eps}(u_{\eps}^{(1)})v \stackrel{(*)}= I'_{\eps}(u_{\eps}^{(1)})(\tilde{\varphi}(\eps \cdot) v) +(1+\norm{v}_2) o(1)} \\
&=& I'_{\eps}(u_{\eps})(\tilde{\varphi}(\eps \cdot) v) - \left(I'_{\eps}(u_{\eps})(\tilde{\varphi}(\eps \cdot) v) - I'_{\eps}(u_{\eps}^{(1)})(\tilde{\varphi}(\eps \cdot) v) \right) +(1+\norm{v}_2) o(1) \\
&=& J'_{\eps}(u_{\eps})(\tilde{\varphi}(\eps \cdot) v) - \left(I'_{\eps}(u_{\eps}) - I'_{\eps}(u_{\eps}^{(1)}) \right)(\tilde{\varphi}(\eps \cdot) v) +(1+\norm{v}_2) o(1).
\end{eqnarray*}
Indeed, we justify $(*)$ as done in Step 3: notice that $u_{\eps}^{(1)}\equiv 0$ outside $\Omega_{h_0}/\eps$ and $1-\tilde{\varphi}(\eps \cdot)\equiv 0$ in $\Omega_{\frac{3}{2}h_0}/\eps$, so in the annulus $(\Omega_{\frac{3}{2} h_0}/\eps) \setminus (\Omega_{h_0}/\eps)$ both $u_{\eps}^{(1)}$ and $1-\tilde{\varphi}(\eps \cdot)$ are zero; notice also that $\complement(\Omega_{\frac{3}{2} h_0}/\eps)$ and $\Omega_{h_0}/\eps$ get far one from the other as $\eps \to 0$. Thus we have
\begin{eqnarray*}
&&\abs{I'_{\eps}(u_{\eps}^{(1)})((1-\tilde{\varphi}(\eps \cdot))v)} =2 \int_{(\Omega_{h_0}/\eps)\times \complement(\Omega_{\frac{3}{2}h_0}/\eps)} \frac{\abs{u_{\eps}^{(1)}(x)} \abs{(1-\tilde{\varphi}(\eps y))v(y)}}{\abs{x-y}^{N+2s}} dx \, dy \\
&&\leq 2\int_{(\Omega_{h_0}/\eps)\times \complement(\Omega_{\frac{3}{2}h_0}/\eps)} \frac{\abs{u_{\eps}(x)} \abs{v(y)}}{\abs{x-y}^{N+2s}}dx \, dy 
\leq (1+\norm{v}_2) o(1) 
\end{eqnarray*}
where in the last passage we argue as for $(I_1)$ in Step 3. Notice that $o(1)$ does not depend on $v$. 
Thus we obtain
\begin{eqnarray*}
\lefteqn{\abs{I'_{\eps}(u_{\eps}^{(1)})v}} \\
 &\leq& \left( \norm{J'_{\eps}(u_{\eps})}_{H^{\text{-}s}(\R^N)} +\norm{I'_{\eps}(u_{\eps}) - I'_{\eps}(u_{\eps}^{(1)})}_{H^{\text{-}s}(\R^N)}\right) \norm{\tilde{\varphi}(\eps \cdot)v}_{H^s(\R^N)} \\
 &&+ (1+\norm{v}_{H^s(\R^N)}) o(1) \\
 &\leq& \left( \norm{J'_{\eps}(u_{\eps})}_{H^{\text{-}s}(\R^N)} + \norm{I'_{\eps}(u_{\eps}) - I'_{\eps}(u_{\eps}^{(1)})}_{H^{\text{-}s}(\R^N)}\right) \norm{v}_{H^s(\R^N)} \big(C +o(1)\big)+\\
 &&+(1+\norm{v}_{H^s(\R^N)}) o(1)
\end{eqnarray*}
where in the last inequality we argue as in Step 5 (again $o(1)$ does not depend on $v$). 
Concluding, we have, by choosing $\norm{v}_{H^s(\R^N)}=1$, that
\begin{eqnarray*}
\lefteqn{\norm{I'_{\eps}(u_{\eps}^{(1)})}_{H^{\text{-}s}(\R^N)}}\\
&& \leq C \left( \norm{J'_{\eps}(u_{\eps})}_{H^{\text{-}s}(\R^N)} + \norm{I'_{\eps}(u_{\eps}) - I'_{\eps}(u_{\eps}^{(1)})}_{H^{\text{-}s}(\R^N)}\right) \big(1 +o(1)\big) + o(1) \to 0
\end{eqnarray*}
by using \eqref{eq_PS2} and \eqref{eq_info_I1}.

\medskip

\noindent \textbf{Step 7.} \textit{Weak convergence of $u_{\eps}^{(1)}$.}

Set $q_{\eps}:=\Upsilon(u_{\eps})$ to avoid cumbersome notation. Since $(\eps q_{\eps})_{\eps} \subset \Omega[0, \nu_0]\subset \Omega$ bounded in $\R^N$, we have that up to a subsequence
$$\eps q_{\eps} \to p_0 \in \overline{\Omega[0, \nu_0]} \subset K_d
\subset \Omega.$$
Moreover, by estimating the norm of $u_{\eps}^{(1)}$ with the norm of $u_{\eps}$ (as done before, in Step 5, for $u_{\eps}^{(2)}$), where $u_{\eps}$ belongs to $S(r_2')$ bounded in $H^s(\R^N)$, we have that also $u_{\eps}^{(1)}$ is a bounded sequence, and thus is so $u_{\eps}^{(1)}(\cdot + q_{\eps})$, which implies, up to a subsequence
$$u_{\eps}^{(1)}(\cdot + q_{\eps}) \wto \tilde{U} \; \; \textnormal{in $H^s(\R^N)$} \quad \textnormal{ as $\eps \to 0$}.$$
For each $v\in H^s(\R^N)$ we apply this weak convergence to the following equalities, derived from \eqref{eq_conv_I1},
\begin{eqnarray*}
\lefteqn{ o(1) = I'_{\eps}(u_{\eps}^{(1)}) v(\cdot-q_{\eps}) = \int_{\R^N} (-\Delta)^{s/2} u_{\eps}^{(1)}(y+q_{\eps}) (-\Delta)^{s/2}v(y) dy +} \\
&&+\int_{\R^N} V(\eps y + \eps q_{\eps}) u_{\eps}^{(1)}(y+q_{\eps}) v(y) dy - \int_{\R^N} f(u_{\eps}^{(1)}(y+q_{\eps})) v(y) dy\\
&=:& (I) + (II) + (III).
\end{eqnarray*}
For $(I)$ and $(III)$ we obtain by the weak convergence and classical arguments 
$$(I)\to \int_{\R^N} (-\Delta)^{s/2} \tilde{U} \,(-\Delta)^{s/2}v \, dy \quad \textnormal{ and } \quad (III)\to - \int_{\R^N} f(\tilde{U}) v \, dy \quad \textnormal{ as $\eps \to 0$}.$$
For $(II)$ instead we have
\begin{eqnarray*}
\lefteqn{\abs{(II) -\int_{\R^N} V(p_0) \tilde{U}v \, dy}}\\ 
&\leq& \Bigg(\int_{\R^N} (V(\eps y + \eps q_{\eps}) -V(p_0))^2 v^2(y) dy\Bigg)^{1/2} \norm{u_{\eps}^{(1)}(\cdot+q_{\eps})}_2+ \\
&&+V(p_0)\Bigg|\int_{\R^N} u_{\eps}^{(1)}(y+q_{\eps}) v(y) dy - \int_{\R^N} \tilde{U}v \, dy \Bigg|
\end{eqnarray*}
where the first term goes to zero (thanks to the boundedness of $u_{\eps}^{(1)}$) by the dominated convergence theorem, while the second thanks to the weak convergence. Thus we finally obtain
$$L'_{V(p_0)}(\tilde{U})v = \int_{\R^N} (-\Delta)^{s/2} \tilde{U} \, (-\Delta)^{s/2}v \, dy+\int_{\R^N} V(p_0) \tilde{U}v \, dy- \int_{\R^N} f(\tilde{U}) v \, dy=0$$
for each $v\in H^s(\R^N)$, that is
\begin{equation}\label{eq_critic_Vp0}
L'_{V(p_0)}(\tilde{U})=0.
\end{equation}

\medskip

\noindent \textbf{Step 8.} \textit{Strong convergence of $u_{\eps}^{(1)}$.}

We want to show the strong convergence of $u_{\eps}^{(1)}(\cdot+q_{\eps})$, that is
\begin{equation}\label{eq_dim_strong_conv_tildeU}
u_{\eps}^{(1)}(\cdot + q_{\eps}) \to \tilde{U} \quad \textnormal{ in $H^s(\R^N)$} \quad \textnormal{ as $\eps \to 0$}.
\end{equation}
Set $\tilde{w}_{\eps}:= u_{\eps}^{(1)}(\cdot+q_{\eps}) -\tilde{U} \wto 0$, again by \eqref{eq_conv_I1} we have
\begin{eqnarray*}
o(1) &=& I'_{\eps}(u_{\eps}^{(1)})\tilde{w}_{\eps}(\cdot-q_{\eps}) \\
&=& L'_{V(p_0)}(\tilde{U})\tilde{w}_{\eps} +\left(\norm{(-\Delta)^{s/2} \tilde{w}_{\eps}}_2^2 + \int_{\R^N} V(\eps y + \eps q_{\eps}) \tilde{w}_{\eps}^2 dy\right)+\\
&&+\int_{\R^N} (V(\eps y + \eps q_{\eps})-V(p_0))\tilde{U}\tilde{w}_{\eps}dy + \int_{\R^N} (f(\tilde{U})-f(\tilde{U}+\tilde{w}_{\eps}))\tilde{w}_{\eps} dy\\
&=:& \Bigg(\norm{(-\Delta)^{s/2} \tilde{w}_{\eps}}_2^2 + \int_{\R^N} V(\eps y + \eps q_{\eps}) \tilde{w}_{\eps}^2 dy\Bigg) + (I) \\
&\geq& \norm{(-\Delta)^{s/2} \tilde{w}_{\eps}}_2^2 +\underline{V} \norm{\tilde{w}_{\eps}}_2^2 + (I)
\end{eqnarray*}
where we have used \eqref{eq_critic_Vp0}. We obtain by the boundedness of $V$ and \eqref{eq_prop_f}
\begin{eqnarray*}
(I) &\geq& - 2 \norm{V}_{\infty} \int_{\R^N}|\tilde{U}| |\tilde{w}_{\eps}|dy-\int_{\R^N} \left(2\beta |\tilde{U|} + C_{\beta}(2^p+1)|\tilde{U}|^p\right)|\tilde{w}_{\eps}| dy -\\
&&- \int_{\R^N} \left(\beta |\tilde{w}_{\eps}|^2 + 2^p C_{\beta} |\tilde{w}_{\eps}|^{p+1}\right) dy
= o(1) - \beta \norm{\tilde{w}_{\eps}}_2^2 - 2^p C_{\beta} \norm{\tilde{w}_{\eps}}_{p+1}^{p+1};
\end{eqnarray*}
in the last passage we have used that $\tilde{w}_{\eps}\wto 0$ in $H^s(\R^N)$, thus by classical properties 
$|\tilde{w}_{\eps}|\wto 0$ in $H^s(\R^N)$ and hence in $L^2(\R^N)$ and in $L^{p+1}(\R^N)$ (observing that $\tilde{U}^p \in L^{1+\frac{1}{p}}(\R^N)$).
Merging together all the things we have, by \eqref{eq_stima_2p} and choosing $\beta < \frac{1}{2} \underline{V}$,
$$o(1) \geq \norm{(-\Delta)^{s/2} \tilde{w}_{\eps}}_2^2 +(\underline{V}-\beta) \norm{\tilde{w}_{\eps}}_2^2 - 2^p C_{\beta} \norm{\tilde{w}_{\eps}}_{p+1}^{p+1}\geq C \norm{\tilde{w}_{\eps}}_{H^s(\R^N)}^2$$
and thus $\tilde{w}_{\eps} \to 0$ strongly in $H^s(\R^N)$, that is the claim.

\medskip

\noindent \textbf{Step 9.} \textit{Localization.}

Observe first that $\tilde{U} \nequiv 0$. Indeed, if not, by \eqref{eq_conv_u2}, \eqref{eq_dim_strong_conv_tildeU} and translation invariance of the norm we would have 
\begin{eqnarray*}
r^*&\leq& \liminf_{\eps \to 0} \norm{U_{\eps}}_{H^s(\R^N)} \\
&\leq& \liminf_{\eps \to 0} \norm{U_{\eps}(\cdot-p_{\eps})+\varphi_{\eps}}_{H^s(\R^N)} + \liminf_{\eps \to 0} \norm{\varphi_{\eps}}_{H^s(\R^N)} \\
&\leq& \lim_{\eps \to 0} \Big( \norm{u_{\eps}^{(1)}}_{H^s(\R^N)} + \norm{u_{\eps}^{(2)}}_{H^s(\R^N)}\Big)+ r_2' \leq r_2' <r^*,
\end{eqnarray*}
impossible. 
By \eqref{eq_dim_strong_conv_tildeU} we obtain also
$$I_{\eps}(u_{\eps}^{(1)}) \to L_{V(p_0)}(\tilde{U}) \quad \textnormal{ as $\eps \to0$}.$$
Thus we find, by using also \eqref{eq_confr_JI} and \eqref{eq_PS1},
$$L_{V(p_0)}(\tilde{U}) = I_{\eps}(u_{\eps}^{(1)}) + o(1) = J_{\eps}(u_{\eps}) + o(1) \leq l_0'+ o(1)$$
and hence, letting $\eps \to 0$,
\begin{equation}\label{eq_stima_Lvp0}
L_{V(p_0)}(\tilde{U}) \leq l_0'< E_{m_0+\nu_1}.
\end{equation}
Moreover by \eqref{eq_critic_Vp0} and $\tilde{U}\nequiv 0$, we have
$$E_{V(p_0)}\leq L_{V(p_0)}(\tilde{U});$$
joining together the two previous inequalities we find $E_{V(p_0)} < E_{m_0+\nu_1}$ which implies, by the monotonicity of $E_a$, that
$$V(p_0) < m_0+\nu_1.$$
Joining this information to the fact that $p_0 \in \Omega$ (and in particular $V(p_0) \geq m_0$) we have $p_0 \in \Omega[0, \nu_1)$, that is
\begin{equation}\label{eq_local_qeps}
\eps \Upsilon(u_{\eps}) \to p_0 \in \Omega[0,\nu_1) \quad \textnormal{ as $\eps \to 0$}.
\end{equation}
Exploiting again \eqref{eq_stima_Lvp0} (observe that $E_{m_0+\nu_1}< E_{m_0+\nu_0} < l_0$) together with $L'_{V(p_0)}(\tilde{U})=0$, and $\tilde{U}\neq 0$, we have that $\tilde{U}$ belongs to $S_{V(p_0)}$ up to translations, that is
$$U:= \tilde{U}(\cdot - y_0) \in S_{V(p_0)}\subset \widehat{S}$$
for some suitable $y_0 \in \R^N$. So, set
$$p_{\eps}:=q_{\eps} + y_0$$
we have
\begin{equation}\label{eq_conv_u1}
\norm{u_{\eps}^{(1)} - U(\cdot - p_{\eps})}_{H^s(\R^N)} \to 0 \quad \textnormal{ as $\eps \to 0$}.
\end{equation}
For a later use observe also that
\begin{equation}\label{eq_conv_p_eps}
\eps p_{\eps} \to p_0 \in \Omega[0, \nu_1) \quad \textnormal{ as $\eps \to 0$}.
\end{equation}

\medskip

\noindent \textbf{Step 10.} \textit{Conclusions.}

By \eqref{eq_local_qeps} we have that
$$\eps \Upsilon(u_{\eps})\in \Omega[0,\nu_1) $$
definitely for $\eps$ small. This is the first part of the claim. 
Moreover, by \eqref{eq_conv_u1} and \eqref{eq_conv_u2} we have
\begin{equation}\label{eq_conv_uu}
\norm{u_{\eps} - U(\cdot - p_{\eps})}_{H^s(\R^N)} \to 0 \quad \textnormal{ as $\eps \to0$}
\end{equation}
and thus, since $\widehat{\rho}(u_{\eps})\leq \norm{u_{\eps} - U(\cdot - p_{\eps})}_{H^s(\R^N)}$ by definition, also $\widehat{\rho}(u_{\eps})\to 0$ and hence
$$\widehat{\rho}(u_{\eps})\in [0, \rho_1]$$
definitely for $\eps$ small. This concludes the proof.
\QED

\bigskip

In the next proposition we see that solutions of $J'_{\eps}(u)=0$ are, under suitable assumptions, also solutions of $I'_{\eps}(u)=0$.

\begin{corollary}\label{corol_passaggio_critici}
Let $(u_{\eps})_{\eps}$ be a sequence of critical points of $J_{\eps}$, that is $J_{\eps}'(u_{\eps})=0$,
satisfying 
$$u_{\eps}\in S(r_2'), \quad J_{\eps}(u_{\eps})\leq l_0' \quad \textit{and} \quad \eps \Upsilon(u_{\eps})\in \Omega[0,\nu_0] $$
for any $\eps>0$. 
Then, for $\eps$ sufficiently small, we have
$$Q_{\eps}(u_{\eps})=0, \quad \textit{and} \quad Q'_{\eps}(u_{\eps})=0.$$
In particular $I'_{\eps}(u_{\eps})=0$, 
which means that $u_{\eps}$ is a solution of \eqref{eq_princ_epsx}.
\end{corollary}

\claim Proof.
By the proof of Lemma \ref{lemma_stima_basso}, we notice, since $1-\varphi_{i_\eps}^{\eps}\equiv 1$ outside $\Omega_{i_{\eps}+1}^{\eps}$, and thus outside $\Omega_{h_0}/{\eps}$, that
\begin{equation}\label{eq_dim_uguagl_norme}
\norm{u_{\eps}}_{L^2(\R^N \setminus (\Omega_{h_0}/\eps))} = \norm{u_{\eps}^{(2)}}_{L^2(\R^N \setminus (\Omega_{h_0}/\eps))} \leq \norm{u_{\eps}^{(2)}}_{H^s(\R^N)}
\end{equation}
and hence $\norm{u_{\eps}}_{L^2(\R^N \setminus (\Omega_{h_0}/\eps))}\to 0$ by \eqref{eq_conv_u2}. 

Through a careful analysis of the Steps 3--5 of the proof, that is by \eqref{eq_dim_stima_mista} and \eqref{eq_dim_uguagl_norme}, we see, more precisely, that
$$\norm{u_{\eps}}^2_{L^2(\R^N \setminus (\Omega_{h_0}/\eps))}\leq \frac{C}{n_{\eps}}+ \frac{C}{n_{\eps}^{2s}} +o(1)$$
where $C=C(r_2')$ and $o(1)$ depends on the rate of convergence of $J'_{\eps}(u_{\eps})$. Thus, since we assume $J'_{\eps}(u_{\eps})\equiv 0$, we gain uniformity, i.e., called $\alpha^*:=\min\{1, 2s\}$, we obtain
\begin{equation}\label{eq_conv_u_a0_unif}
\norm{u_{\eps}}^2_{L^2(\R^N \setminus (\Omega_{h_0}/\eps))} \leq \frac{C}{n_{\eps}^{\alpha^*}} \sim \eps^{\alpha^*/2}
\end{equation}
As a consequence
$$\frac{1}{\eps^{\alpha}}\norm{u_{\eps}}^2_{L^2(\R^N \setminus (\Omega_{h_0}/\eps))} \to 0 \quad \textnormal{ as $\eps \to 0$}$$
for $\alpha \in (0, \alpha^*/2)$, and hence $Q_{\eps}(u_{\eps})\equiv Q'_{\eps}(u_{\eps})\equiv 0$ for $\eps$ sufficiently small.
\QED

\bigskip

We want to show now a (truncated) Palais-Smale-like condition.

\begin{proposition}\label{prop_Palais_Smale}
There exists $r_2'' \in (0, \min\{r_0, r_1\})$ sufficiently small with the following property: let $\eps>0$ fixed and let $(u_j)_{j}\subset S(r_2'')$ be such that 
\begin{equation}\label{eq_PS3}
\norm{J'_{\eps}(u_j)}_{H^{\text{-}s}(\R^N)}\to 0 \quad \textit{ as $j\to +\infty$}
\end{equation}
with the additional assumption
$$(\eps\Upsilon(u_j))_j \subset \Omega[0,\nu_0].$$
Then $(u_j)_j$ admits a strongly convergent subsequence in $H^s(\R^N)$.
\end{proposition}

\medskip

\claim Proof.
Let $r_2''$ to be fixed. Since $S(r_2'')$ is bounded, up to a subsequence we can assume $u_j \wto u_0$ in $H^s(\R^N)$. We want to show that
$$\lim_{R \to +\infty} \lim_{j \to +\infty} \norm{u_j}_{L^q(\R^N \setminus B_R)}=0$$
for $q=2$ and $q=p+1$ and conclude by Lemma \ref{lemma_conv_forte}. 

Arguing similarly to Step 1 of the proof of Lemma \ref{lemma_stima_basso}, i.e. exploiting Remark \ref{remark_triplanorma}, we obtain for $L \gg 0$, uniformly in $j \in \N$,
$$\tnorm{u_j}_{\R^N \setminus B_L} \leq C r_2'';$$
indeed we work with the set $B_L - \Upsilon(u_j)$ which expands to $\R^N$ as $L, j\to +\infty$, since $\Upsilon(u_j) \in \Omega/\eps$, a fixed bounded set. 
Moreover, for any $n \in \N$, we have
$$\sum_{i=1}^n \norm{u_j}^2_{L^2(B_{L+ni}\setminus B_{L+n(i-1)})} \leq (C r_2'')^2$$
and similarly for the Gagliardo seminorm and the $(p+1)$-norm, thus for some $i_{j,n}\in \{1, \dots, n\}$
$$ \norm{u_j}^2_{ A^{j,n}} + \norm{u_j}^{p+1}_{L^{p+1}(A^{j,n})} \leq \frac{C}{n}$$
where $A^{j,n}:=B_{L+ni_{j,n}}\setminus B_{L+n(i_{j,n}-1)}$. Again similarly to Step 2 of the proof of Lemma \ref{lemma_stima_basso}, we introduce $\psi_{j,n}$ such that
$$B_{L + n(i_{j,n} -1)} \prec \psi_{j,n} \prec B_{L+ni_{j,n}}$$ 
and $\norm{\nabla \psi_{j,n}}_{\infty}=o(1)$ as $n\to +\infty$;
moreover we set
$$\tilde{u}_{j,n} :=(1- \psi_{j,n}) u_j.$$
Observe that $\chi_{B_L} \leq \psi_{j,n}$, thus $\supp(\tilde{u}_{j,n}) \subset \complement(B_L)$. 
Arguing as in Step 5 and 3 of the proof of Lemma \ref{lemma_stima_basso} we obtain
\begin{eqnarray*}
\lefteqn{ 
\int_{\R^{2N}} \frac{\abs{\tilde{u}_{j,n}(x)-\tilde{u}_{j,n}(y)}^2}{|x-y|^{N+2s}} \leq 4 \int_{\complement(B_L) \times \R^{N}} \frac{\abs{\psi_{j,n}(x)-\psi_{j,n}(y)}^2\abs{u_j(x)}^2}{|x-y|^{N+2s}} dx \, dy +}\\
 &&+4 \int_{\complement(B_L)\times \R^N} \frac{\abs{u_j(x)-u_j(y)}^2}{|x-y|^{N+2s}} dx \, dy \leq o(1) \norm{u_j}_{L^2(\complement(B_L))}^2 + C [u_j]_{\complement(B_L), \R^N}^2 
\end{eqnarray*}
thus $\norm{\tilde{u}_{j,n}}_{H^s(\R^N)} \leq C r_2''$ and hence, choosing $r_2''$ sufficiently small, we have 
$$\norm{\tilde{u}_{j,n}}_{H^s(\R^N)} \leq r_1;$$
by Lemma \ref{lemma_stima_I}, for $q \in \{2, p+1\}$, we obtain
$$\norm{u_j}_{L^q(\R^N \setminus B_{L+ni_{j,n}})}^2 = \norm{\tilde{u}_{j,n}}_{L^q(\R^N \setminus B_{L+ni_{j,n}})}^2 \leq C \norm{\tilde{u}_{j,n}}_{H^s(\R^N)}^2 \leq C I'_{\eps}(\tilde{u}_{j,n})\tilde{u}_{j,n}.$$
Thus the claim comes if we show that 
$$I'_{\eps}(\tilde{u}_{j,n})\tilde{u}_{j,n} \to 0 \quad \textnormal{ as $j,n\to +\infty$}.$$
Indeed we have
\begin{eqnarray*}
\lefteqn{I'_{\eps}(\tilde{u}_{j,n})\tilde{u}_{j,n} = I'_{\eps}(u_j) \tilde{u}_{j,n} - \int_{\R^{2N}} (-\Delta)^{s/2}(\psi_{j,n} u_j) (-\Delta)^{s/2}((1-\psi_{j,n})u_j) dx-}\\
 && -\int_{A^{j,n}} V(\eps x) \psi_{j,n} (1-\psi_{j,n}) u_j^2dx - \int_{A^{j,n}} (f((1-\psi_{j,n})u_j)- f(u_j)) (1-\psi_{j,n})u_j dx\\
&&=: J'_{\eps}(u_j)\tilde{u}_{j,n} - Q'_{\eps}(u_j)\tilde{u}_{j,n} + (I)
\leq o(1) + (I)
\end{eqnarray*}
where we have used that $J'_{\eps}(u_j) \to 0$ (as $j\to +\infty$, uniformly in $n\in \N$), the boundedness of $\norm{\tilde{u}_{j,n}}_{H^s(\R^N)}$ and the positivity of $ Q'_{\eps}(u_j)\tilde{u}_{j,n}$. The term $(I)$ can be estimated in the same way as done in Steps 3-4 of the proof of Lemma \ref{lemma_stima_basso} (fixed $j\in \N$, and $n\to +\infty$), and hence we reach the claim.
\QED

\subsection{Deformation lemma on a neighborhood of expected solutions}

We want to define now a neighborhood of expected solutions (see \cite{CJT}), which will be invariant under a suitable deformation flow. 
Consider $r_3:= \min\{ r', r_0', r_2', r_2''\}$ (see Lemma \ref{lemma_stima_g}, Lemma \ref{lemma_stimabasso_I}, Theorem \ref{theorem_stima_J'} and Proposition \ref{prop_Palais_Smale}), and let us define
$$R(\delta, u):= \delta - \frac{\delta_2}{2} (\widehat{\rho}(u) - \rho_1)_+ \leq \delta$$
and
$$\mathcal{X}_{\eps, \delta}:= \Big\{ u \in S(\rho_0) \mid \eps \Upsilon(u) \in \Omega[0,\nu_0), \; J_{\eps}(u) < E_{m_0} + R(\delta, u) \Big\}$$
where
$$0< \rho_1 < \rho_0 < r_3,$$
$\eps$ is sufficiently small and 
\begin{equation}\label{eq_propr_delta}
\delta \in \left(0, \min \left\{ \tfrac{\delta_2}{4}(\rho_0- \rho_1), \, \delta_1, \, l_0'-E_{m_0}\right\}\right);
\end{equation}
here $\delta_1$ and $\delta_2$ are the ones that appear in Lemma \ref{lemma_stimabasso_I} and Theorem \ref{theorem_stima_J'}. Notice that the \emph{height} of the sublevel in $\mc{X}_{\eps, \delta}$ depends on $u$ itself; this will be used to gain a deformation which preserves $\mc{X}_{\eps, \delta}$. 

\noindent We begin by pointing out some geometrical features of the neighborhood $\mc{X}_{\eps, \delta}$. 
\begin{itemize}
\item $\mc{X}_{\eps, \delta}$ is open. Indeed, $S(\rho)$ and $\{J_{\eps}(u) < E_{m_0} + R(\delta, u)\}$ are open, and $\Omega[0,\nu_0)= \Omega(-\gamma, \nu_0)$ for a whatever $\gamma>0$ (since $V$ cannot go under $m_0$ in $\Omega$) and thus open. Moreover it is nonempty (see e.g. Section \ref{sec_def_phi_eps}).
\item If $v\in\mathcal{X}_{\eps,\delta} \subset S(\rho_0)$, then by \eqref{eq_stima_widerho} we have 
$\widehat{\rho}(v)<\rho_0.$

\item 
If $v\in \mathcal{X}_{\eps, \delta} \subset \{ \eps \Upsilon(v) \in \Omega[0, \nu_0]\}$, then
\begin{equation}\label{eq_appart_insieme_piccolo}
\eps \Upsilon(v) \in \Omega [0, \nu_1).
\end{equation}
Indeed, if not, i.e. $\eps \Upsilon(v) \in \Omega [\nu_1, \nu_0]$, then by Lemma \ref{lemma_stimabasso_I} we have
$$J_{\eps}(v) \geq I_{\eps}(v) \geq E_{m_0} + \delta_1 > E_{m_0} + \delta \geq E_{m_0} + R(\delta, v)$$
which is an absurd.

\item If $R(\delta, v) \geq -\delta$ then
\begin{equation}\label{eq_cons_conten}
v \in S(\rho_0).
\end{equation}
Indeed $\frac{\delta_2}{4}(\widehat{\rho}(u)-\rho_1)_+ \leq \delta$ implies, by the restriction on $\delta$,
$$(\widehat{\rho}(u)-\rho_1)_+ < \rho_0-\rho_1.$$
If $\widehat{\rho}(u)<\rho_1$ then clearly $u \in S(\rho_1)\subset S(\rho_0)$. If instead $\widehat{\rho}(u)\geq \rho_1$, then 
$\widehat{\rho}(u) < \rho_0$, 
which again implies $u \in S(\rho_0).$
\end{itemize} 
We further define the set of critical points of $J_{\eps}$ lying in the neighborhood of expected solutions
$$K_c:= \left\{ u \in \mathcal{X}_{\eps, \delta} \mid J'_{\eps}(u)=0, \; J_{\eps}(u)=c\right\},$$
the sublevel
 $$\mc{X}_{\eps, \delta}^c:= \mc{X}_{\eps, \delta}\cap J_{\eps}^c$$
 and
 $$ (\mc{X}_{\eps, \delta})^c_d:= \{ u \in \mc{X}_{\eps, \delta} \mid d \leq J_{\eps}(u) \leq c\}, $$
for every $c,d \in \R$. 
We present now a deformation lemma with respect to $K_c$, for $c$ sufficiently close to $E_{m_0}$. 

\begin{lemma}\label{lem_def_lem}
Let $c \in (E_{m_0}-\delta, E_{m_0}+\delta)$.
Then there exists a deformation at level $c$, which leaves the set $\mathcal{X}_{\eps, \delta}$ invariant.
That is, for every $U$ neighborhood of $K_c$ ($U=\emptyset$ if $K_c=\emptyset$), there exist a small $\omega>0$ and a continuous deformation $\eta: [0,1]\times \mc{X}_{\eps, \delta} \to \mc{X}_{\eps,\delta}$ such that
\begin{itemize}
\item[(i)] \ $\eta(0, \cdot)=id$;
\item[(ii)] \, $J_{\eps}(\eta(\cdot, u))$ is non-increasing;
\item[(iii)] \, $\eta(t, u)=u$ for every $t \in [0,1]$, if $J_{\eps}(u) \notin (E_{m_0}-\delta, E_{m_0}+\delta)$;
\item[(iv)] \, $\eta(1,\mc{X}_{\eps, \delta}^{c+\omega}\setminus U)\subset \mc{X}_{\eps, \delta}^{c-\omega}$;
\item[(v)] \, $\eta(\cdot, u)$ is a semigroup.
\end{itemize}
\end{lemma}

\claim Proof. 
Let $\mathcal{V}: \{ u \in H^s(\R^N) \mid J'_{\eps}(u)\neq 0\} \to H^s(\R^N)$ be a locally Lipschitz pseudo-gradient vector field associated to $J_{\eps}$, and let $\phi \in Lip_{loc}(H^s(\R^N), \R)$ be a cutoff function such that 
$\supp (\phi) \subset (\mc{X}_{\eps, \delta})^{E_{m_0}+\delta}_{E_{m_0}-\delta}$ 
and $\phi=1$ in a small neighborhood of $c$. 
We consider the Cauchy problem
\begin{equation}\label{eq_dim_CP}
\parag{&\dot{\eta} = - \phi(\eta) \frac{\mathcal{V}(\eta)}{\norm{\mathcal{V}(\eta)}_{H^s(\R^N)}},& \\ &\eta(0, u)=u.&}
\end{equation}
The proof keeps on classically, obtaining a deformation $\eta: [0,1]\times \mathcal{X}_{\eps, \delta} \to H^s(\R^N)$. 
We want to prove now that $\eta$ goes into $\mathcal{X}_{\eps, \delta}$.

Let $u\in \mathcal{X}_{\eps, \delta}$. We need to show that $\eta(t,u)\in \mathcal{X}_{\eps, \delta}$ for every $t>0$. 
Since $\mathcal{X}_{\eps, \delta}$ is open, $\eta(s,u)$ continues staying in $\mathcal{X}_{\eps, \delta}$ for $s$ small. Thus assume that 
$$\eta(s,u) \in \mathcal{X}_{\eps,\delta}, \quad \textnormal{ for every $0\leq s<t_0$}$$
for some $t_0>0$, and we want to show that $\eta(t_0, u)\in\mathcal{X}_{\eps, \delta}$. 
Notice first that, by using (iii), (v), (ii) and the continuity of $\eta$ and $J_{\eps}$ 
we can assume that 
\begin{equation}\label{eq_dim_stima[)}
J_{\eps}(\eta(t_0, u)) \in [E_{m_0}-\delta, E_{m_0}+\delta).
\end{equation}

\textbf{Step 1: $\eps \Upsilon(\eta(t_0,u)) \in \Omega[0, \nu_0)$.} By \eqref{eq_appart_insieme_piccolo} we have
$$\eps \Upsilon(\eta(s,u)) \in \Omega[0,\nu_1), \quad \textnormal{ for every $0\leq s<t_0$},$$
and thus by continuity 
$\eps \Upsilon(\eta(t_0,u)) \in \overline{\Omega[0,\nu_1)} \subset \Omega[0,\nu_0].$

\textbf{Step 2: $J_{\eps}(\eta(t_0,u)) < E_{m_0} + R(\delta, \eta(t_0, u))$.} 
If $\widehat{\rho}(\eta(t_0, u)) \leq \rho_1$ then $R(\delta, \eta(t_0, u))=\delta$ and we directly have the claim, recalled that $J_{\eps}(\eta(t_0, u))<E_{m_0}+\delta$ by \eqref{eq_dim_stima[)}. 
Assume instead $\widehat{\rho}(\eta(t_0, u)) > \rho_1$. By continuity, there exists $t_1 \in (0, t_0)$ such that we have
$$\widehat{\rho}(\eta(s, u)) > \rho_1, \quad \textnormal{ for every $s\in [t_1, t_0]$.} $$
In particular 
$$\parag{&\eta(s, u) \in S(\rho_0) \subset S(r_3) \subset S(r_2'),& \\ 
&J_{\eps}(\eta(s, u)) < E_{m_0}+\delta < E_{m_0+\nu_1},& \\ 
&\widehat{\rho}(\eta(s, u)) \in (\rho_1, \rho_0],& \\ 
&\eps \Upsilon(\eta(s, u)) \in \Omega[0, \nu_0]&
}
$$
for $s\in [t_1, t_0]$. 
Then by Theorem \ref{theorem_stima_J'} we have
$$\norm{J'_{\eps}(\eta(s, u))}_{H^{\text{-}s}(\R^N)} \geq \delta_2, \quad \textnormal{ for every $s\in [t_1, t_0]$.} $$
We can thus compute with standard argument, by using \eqref{eq_dim_CP}, the properties of the pseudo-gradient and \eqref{eq_lipsc_rho}, 
\begin{eqnarray*}
J_{\eps}(\eta(t_0, u)) &\leq & J_{\eps}(\eta(t_1, u)) - \frac{\delta_2}{2}\big(\widehat{\rho}(\eta(t_0, u)) - \widehat{\rho}(\eta(t_1, u))\big)\\
&<& E_{m_0} + \delta - \frac{\delta_2}{2}\big(\widehat{\rho}(\eta(t_1,u)) - \rho_1\big) - \frac{\delta_2}{2}\big(\widehat{\rho}(\eta(t_0, u)) - \widehat{\rho}(\eta(t_1, u))\big)\\
&= &E_{m_0} + R(\delta, \eta(t_0, u)),
\end{eqnarray*}
that is the claim.

\textbf{Step 3: $\eta(t_0, u)\in S(\rho_0)$.} 
By the previous point we have
$J_{\eps}(\eta(t_0,u)) \leq E_{m_0} + R(\delta, \eta(t_0, u)).$
Since \eqref{eq_dim_stima[)} implies $J_{\eps}(\eta(t_0,u))\geq E_{m_0}-\delta$, then by \eqref{eq_cons_conten} we have $\eta(t_0,u) \in S(\rho_0)$, 
and thus the claim. 
\QED

\subsection{Maps homotopic to the embedding}

We search now for two maps $\Phi_{\eps}$, $\Psi_{\eps}$ such that, for a sufficiently small $\sigma_0\in (0,1)$ and a sufficiently small $\hat{\delta}=\hat{\delta}(\sigma_0)\in (0, \delta)$ (see \eqref{eq_propr_delta}), defined
$$ I:=[1-\sigma_0, 1+\sigma_0],$$
we have, for small $\eps$,
$$I \times K \; \stackrel{\Phi_{\eps}} \to \; \mathcal{X}_{\eps, \delta}^{E_{m_0}+\hat{\delta}} \;\stackrel{\Psi_{\eps}} \to \; I \times K_d $$
with the additional condition
$$\partial I \times K \;\stackrel{\Phi_{\eps}} \to \;\mathcal{X}_{\eps, \delta}^{E_{m_0}-\hat{\delta}} \;\stackrel{\Psi_{\eps}} \to\; (I\setminus \{1\}) \times K_d;$$
then we will prove that $\Psi_{\eps} \circ \Phi_{\eps}$ is homotopic to the identity. While the first property is useful for category arguments to gain multiplicity of solutions, the second additional condition will be essential for developing \emph{relative} category (and cup-length) arguments.

\subsubsection{Definition of $\Phi_{\eps}$}\label{subsubsec_Phi}
\label{sec_def_phi_eps}

Let us fix a ground state $U_0 \in S_{m_0}\subset \widehat{S}$, i.e. $L_{m_0}(U_0)=E_{m_0}$ (see Theorem \ref{thm_esist_Poh_min} and \eqref{eq_ug_poh_gr.st}). Define, for $p \in K$ and $t \in I$ ($\sigma_0$ to be fixed)
$$\Phi_{\eps}(t,p):= U_0\left(\tfrac{\cdot-p/\eps}{t}\right)\in H^s(\R^N).$$
We show now that, for $\eps$ small, $\Phi_{\eps}(t,p) \in \mathcal{X}_{\eps, \delta}^{E_{m_0}+\hat{\delta}} $.
\begin{itemize}
\item 
$\Phi_{\eps}(t,p) \in S(\rho_1) \subset S(\rho_0)$:

indeed, recalled that the dilation $t\in \R \mapsto U_0(\cdot/t)\in H^s(\R^N)$ is continuous, we have
$$\norm{U_0\left(\tfrac{\cdot-p/\eps}{t}\right)-U_0(\cdot-p/\eps)}_{H^s(\R^N)}= \norm{U_0\left( \cdot/t\right)-U_0}_{H^s(\R^N)}< \rho_1$$
for $t \in I$ and sufficiently small $\sigma_0=\sigma_0(U_0)$ (not depending on $\eps$). Thus, setting $\varphi_t:= U_0\left(\frac{\cdot-p/\eps}{t}\right)-U_0(\cdot-p/\eps)$ we have
$$\Phi_{\eps}(t,p) = U_0(\cdot-p/\eps) + \varphi_t$$
with $U_0 \in \widehat{S}$, $p/\eps \in \R^N$ and $\norm{\varphi_t}_{H^s(\R^N)}<\rho_1$, which is the claim.

\item $\eps \Upsilon(\Phi_{\eps}(t,p)) \in \Omega[0, \nu_0)$:

indeed, by the previous point and Lemma \ref{lem_def_bar}, we have
$$\abs{\Upsilon(\Phi_{\eps}(t,p)) - p/\eps} < 2R_0$$
hence $\abs{\eps \Upsilon(\Phi_{\eps}(t,p)) - p} < 2\eps R_0$, and since $p \in K$
$$d(\eps \Upsilon(\Phi_{\eps}(t,p)), K) < 2\eps R_0.$$
For sufficiently small $\eps$, we have $K_{2\eps R_0} \subset \Omega[0, \nu_0)$, and thus the claim. 
In particular, for a later use observe that
\begin{equation}\label{eq_posizione_Ups}
\eps\Upsilon(\Phi_{\eps}(t,p)) = p + o(1).
\end{equation}

\item $J_{\eps}(\Phi_{\eps}(t,p)) < E_{m_0}+R(\delta, \Phi_{\eps}(t,p))$:

indeed $\Phi_{\eps}(t,p) \in S(\rho_1)$, thus $\widehat{\rho}(\Phi_{\eps}(t,p)) < \rho_1$, which implies $R(\delta, \Phi_{\eps}(t,p) ) = \delta$ and the claim comes from the following point, since $\hat{\delta}<\delta$.

\item $J_{\eps}(\Phi_{\eps}(t,p)) < E_{m_0}+\hat{\delta}$:

indeed we have by Lemma \ref{lemma_stima_g0} (b)
\begin{eqnarray}
\lefteqn{J_{\eps}(\Phi_{\eps}(t,p))\notag}\\ 
&=& L_{m_0}\left(\Phi_{\eps}(t,p)\right) + \frac{1}{2}\int_{\R^N} (V(\eps x)- m_0) \Phi_{\eps}^2(t,p)dx+ Q_{\eps}\left(\Phi_{\eps}(t,p)\right) \notag\\
&=:& L_{m_0}\left(U_0\left(\tfrac{\cdot-p/\eps}{t}\right)\right) + (I) +(II) = g(t) E_{m_0} + o(1)
 \label{eq_dim_stima_g(t)}\\
&\leq & E_{m_0} +o(1) 
\notag
\end{eqnarray}
where we used $g(t) \leq 1$. Indeed, as regards $(I)$ we have
$$(I)= \frac{1}{2}\int_{\R^N} (V(\eps x + p)- m_0) U_0^2(x/t) dx\to 0 \quad \textnormal{ as $\eps \to 0$}$$
by exploiting that $p \in K$ and the dominated convergence theorem, together with the boundedness of $V$. 
Focusing on $(II)$ instead, we have 
$$(II) = \left( \frac{1}{\eps^{\alpha}} \norm{U_0(\cdot/t)}_{L^2(\R^N \setminus ((\Omega_{2h_0} - p)/\eps)}^2 -1 \right)_+^{\frac{p+1}{2}};$$
since $p \in K \subset \Omega \subset \Omega_{2h_0}$, we have $0 \in \Omega_{2h_0}-p$ and moreover $B_r \subset \Omega_{2h_0}-p$ for some ball $B_r$; notice that $B_r/\eps$ covers the whole $\R^N$ as $\eps \to 0$. Therefore, by the polynomial estimate we have
$$ \norm{U_0(\cdot/t)}_{L^2(\R^N \setminus ((\Omega_{2h_0} - p)/\eps)}^2 \leq C \norm{\tfrac{1}{1+|x|^{N+2s}}}_{L^2(\R^N \setminus (B_r/\eps))} ^2 \leq C \eps^{N+4s},$$
and hence $(II) \to 0$ as $\eps \to 0$, since $\alpha < N+4s$. 
Therefore, by choosing a sufficiently small $\eps$, we obtain
$$J_{\eps}(\Phi_{\eps}(t,p)) \leq E_{m_0} + \tfrac{1}{2} \hat{\delta} < E_{m_0} + \hat{\delta}.$$
\end{itemize}
\vspace{-\topsep}
Finally, we show the additional condition.

\begin{itemize}
\item $J_{\eps}(\Phi_{\eps}(1\pm \sigma_0,p)) < E_{m_0}-\hat{\delta}$:

indeed, looking at \eqref{eq_dim_stima_g(t)} we see that, for small $\eps$,
$$J_{\eps}(\Phi_{\eps}(1\pm \sigma_0 ,p)) < g(1\pm \sigma_0) E_{m_0} + \hat{\delta};$$
since $ g(1\pm \sigma_0)<1$, we can find a small $\hat{\delta}< \frac{1-g(1\pm \sigma_0)}{2} E_{m_0}$ (not depending on $\eps$) such that
\begin{equation}\label{eq_buona_post_Phieps}
J_{\eps}(\Phi_{\eps}(1\pm \sigma_0 ,p)) < g(1\pm \sigma_0) E_{m_0} + \hat{\delta} < E_{m_0} - \hat{\delta}
\end{equation}
and thus the claim.
\end{itemize}
\vspace{-\topsep}

\subsubsection{Definition of $\Psi_{\eps}$}

Define a truncation
$$T(t):= \parag{1-\sigma_0 & \quad \textnormal{if $t \leq 1-\sigma_0$},\\ t & \quad \textnormal{if $t \in (1-\sigma_0 ,1+\sigma_0)$}, \\ 1+\sigma_0 & \quad \textnormal{if $t \geq 1+\sigma_0$}}$$
for $t \in \R$, and
$$\Psi_{\eps}(u):= \big( T(P_{m_0}(u)), \eps \Upsilon(u)\big)$$
for every $u \in \mathcal{X}_{\eps, \delta}^{E_{m_0}+\hat{\delta}}$. By the definition of $T$ and property \eqref{eq_appart_insieme_piccolo}, we have directly
$$\Psi_{\eps}(u) \in I \times \Omega[0, \nu_1] \subset I \times \Omega[0, \nu_0] \subset I \times K_d.$$
Assume now $u\in \mathcal{X}_{\eps, \delta}^{E_{m_0}-\hat{\delta}}$. We have, by using Lemma \ref{lemma_stima_J} and Lemma \ref{lemma_stima_g}, 
$$E_{m_0}-\hat{\delta} \geq J_{\eps}(u) \geq L_{m_0}(u) - C_{min} \eps^{\alpha} \geq g(P_{m_0}(u)) E_{m_0} - C_{min}\eps^{\alpha}$$
and hence
$$E_{m_0} \geq g(P_{m_0}(u)) E_{m_0} + \hat{\delta} - C_{min}\eps^{\alpha}> g(P_{m_0}(u)) E_{m_0}$$
where the last inequality holds for $\eps$ small, not depending on $u$. Thus
$$g(P_{m_0}(u))< 1$$
and this must imply, by the properties of $g$, that $P_{m_0}(u) \neq 1$, and in particular
$$T(P_{m_0}(u))\neq 1.$$
This reaches the goal.

\subsubsection{An homotopy to the identity}

Introduce a notation from the algebraic topology: we write, for $ B \subset A$ and $B' \subset A'$,
$$f: (A,B) \to (A', B')$$
whenever
$$f\in C(A,A')\quad \textnormal{ and } \quad f(B) \subset B'.$$
Observed that $\Phi_{\eps}$ and $\Psi_{\eps}$ are continuous, we can rewrite the stated properties as
$$\Phi_{\eps} : \Big( I\times K, \, \partial I \times K \Big) \to \Big( \mathcal{X}_{\eps, \delta}^{E_{m_0}+\hat{\delta}}, \, \mathcal{X}_{\eps, \delta}^{E_{m_0}-\hat{\delta}}\Big),$$
$$\Psi_{\eps}: \Big( \mathcal{X}_{\eps, \delta}^{E_{m_0}+\hat{\delta}}, \, \mathcal{X}_{\eps, \delta}^{E_{m_0}-\hat{\delta}}\Big) \to \Big( I \times K_d, \, (I\setminus \{1\}) \times K_d\Big)$$
and
$$\Psi_{\eps}\circ \Phi_{\eps}: \Big( I \times K, \, \partial I \times K \Big) \to \Big( I \times K_d, \, (I \setminus \{1\}) \times K_d \Big),$$
where a straightforward computation shows
$$ (\Psi_{\eps}\circ \Phi_{\eps}) (t, p) = \left( t, \, \eps \Upsilon\left( U_0\left(\tfrac{\cdot-p/\eps}{t}\right)\right)\right).$$
Clearly, we notice that the inclusion map has the same property, that is set $j(t,p):=(t,p)$ we have
$$j: \Big( I \times K, \, \partial I \times K \Big) \to \Big( I \times K_d, \, (I \setminus \{1\}) \times K_d\Big).$$
We want to show that these maps are homotopic, information useful in the theory of relative cup-length.

\begin{proposition} \label{prop_esist_homot}
For sufficiently small $\eps$, the maps $\Psi_{\eps}\circ \Phi_{\eps}$ and $j$ are homotopic, that is there exists a continuous map 
$H: [0,1] \times I \times K \to I \times K_d$
such that
$$H(\theta, \cdot, \cdot): \Big( I \times K, \, \partial I \times K \Big) \to \Big( I \times K_d, \, (I \setminus \{1\}) \times K_d\Big)$$
for each $\theta \in [0,1]$, with
$H(0, \cdot, \cdot) = \Psi_{\eps}\circ \Phi_{\eps}$ and $ H(1, \cdot, \cdot) = j.$
\end{proposition}

\claim Proof.
Noticed that also $\Psi_{\eps}\circ \Phi_{\eps}$ fixes the first variable, it is sufficient to link the second variables through a segment, that is
$$H(\theta, t, p) := \left( t, \, (1-\theta) \eps \Upsilon\left( U_0\left(\tfrac{\cdot-p/\eps}{t}\right)\right) + \theta p\right),$$
with $\theta \in [0,1]$. We must check that $H$ is well defined, since $K_d$ is not a convex set, generally. Indeed we have, by \eqref{eq_posizione_Ups}
$$(1-\theta) \eps \Upsilon\left( U_0\left(\tfrac{\cdot-p/\eps}{t}\right)\right) + \theta p = (1-\theta) p + o(1) + \theta p = p + o(1).$$
Since $p \in K$, for sufficiently small $\eps$ we have that $p+o(1) \in K_d$, and thus the claim.
\QED

\section{Existence of multiple solutions}\label{sec_cuplenght}

We finally come up to the existence of multiple solutions. 
Here the algebraic notions of \emph{relative category} and \emph{relative cup-length} are of key importance. 

\begin{remark}\label{rem_sola_esist}
Before coming up to multiplicity results, we highlight that existence of a single solution could be obtained without any use of algebraic tools. 
We omit the details, gaining instead the existence as a byproduct of the multiplicity.
\end{remark}

In what follows we refer to \cite{Szu, FLRW} and references therein, and to \cite{Bar} for a complete treatment on the topic of category.

\begin{definition}
Let $X$ be a topological space and let $A, B$ be two closed subsets of $X$. We call the \emph{category of $A$ in $X$, relative to $B$}, and write $k=\cat_{X, B}(A)$, the least integer $k\in \N$ such that there exist $A_0, A_1, \dots, A_k$ closed subsets of $X$ which verify
\begin{itemize}
\item $(A_i)_{i=0\dots k}$ cover $A$;
\item $(A_i)_{i=1 \dots k}$ are contractible;
\item $A_0$ is deformable in $B$, i.e. there exists a continuous $h_0:[0,1]\times (A_0 \cup B) \to X$ such that $h_0(0, \cdot)=id$, $h_1(1, A_0)\subset B$ and $h_0(t, B)\subset B$ for each $t \in [0,1]$.
\end{itemize}
If such $k$ does not exists, we set $\cat_{X, B}(A):=+\infty$. Notice that $\cat_X(A)=\cat_{X, \emptyset}(A)$.
\end{definition}

Since we will work only with $X=\mc{X}_{\eps, \delta}$, to avoid cumbersome notation we will write $\cat(A):= \cat_{\mc{X}_{\eps, \delta}}(A)$ and $\cat(A,B):=\cat_{\mc{X}_{\eps, \delta}, B}(A)$.

We recall now the definition of \emph{(relative) cup-length} (see e.g. \cite{Cha} and references therein). Let $B\subset A$ be topological spaces, and let $\mathbb{F}$ be a whatever field. 
Let moreover $H^*(A)=\bigoplus_{q\geq 0} H^q(A)$ be the Alexander-Spanier cohomology with coefficients in some field $\mathbb{F}$ and let $H^*(A, B)$ be the \emph{relative} Alexander-Spanier cohomology (see \cite{Mas} and references therein). Let $\smile: H^*(A)\times H^*(A)\to H^*(A)$ be the cup-product and let $\smile: H^*(A, B) \times H^*(A) \to H^*(A, B)$ be the relative cup-product. 
In addition recall that a whatever continuous function $f:(A,B) \to (A',B')$, being $(A', B')$ another topological pair, induces homomorphisms $f^*: H^*(A) \to H^*(A)$ and $f^*: H^*(A', B') \to H^*(A, B)$. 
We define
\begin{eqnarray*}
\cupl(f):= \max \{ l \in \N \; | &\exists \alpha_i \in H^{q_i}(A'), \, q_i \geq 1, \textnormal{ for $i=1 \dots l$,} \; \exists \alpha_0 \in H^*(A', B'),&\\
& \textnormal{s.t.} \; f^*(\alpha_0) \smile f^*(\alpha_1) \smile \dots \smile f^*(\alpha_l) \neq 0 \, \textnormal{ in $H^*(A,B)$} \};&
\end{eqnarray*}
if such $l\in \N$ does not exist, but $f^*\nequiv 0$, it results $\cupl(f):=0$, otherwise we define $\cupl(f):=-1$. Moreover, we set $ \cupl(A, B) :=\cupl(id_{(A, B)})$ and $\cupl(A):=\cupl(A, \emptyset)$. 
In the case $A$ is not connected, a slightly different definition (which makes the cup-length additive) can be found in \cite{Bar}. 
Here we collect some properties of the cup-length, see \cite[Lemma 2.6]{BW}.

\begin{lemma}\label{lem_propr_cupl}
We have the following properties.
\begin{itemize}
\item[(a)] For any $f:(A,B)\to (A', B')$ and $f':(A', B')\to (A'', B'')$ it results that
$$\cupl(f' \circ f) \leq \min\{ \cupl(f'), \cupl (f)\}.$$
As a consequence,
\begin{equation}\label{eq_propr_cupl}
\cupl(f' \circ f) \leq \cupl(A', B').
\end{equation}
\item[(b)] For any $f, g:(A,B) \to (A', B')$ homotopic, we have
$$\cupl(f)=\cupl(g).$$
\end{itemize}
\end{lemma}

Finally, we cite the following result \cite{BarP} which can be found in \cite[Lemma 5.5]{CJT}.
\begin{lemma}\label{lem_cupl_jK}
Consider the inclusion $j:(I\times K, \partial I \times K) \to (I\times K_d, \partial I \times K_d)$ for a whatever $K\subset \R^N$ compact, $K_d=\{ x \in \R^N \mid d(x,K)\leq d\}$, and $I=[a,b]$. Then, for $d>0$ sufficiently small, we have
$$\cupl(j) \geq \cupl(K).$$
\end{lemma}

The following lemma links the concepts of category and cup-length, and it can be found in \cite[Proposition 2.6]{Szu} or \cite[Theorem 3.6]{FLRW}.
\begin{lemma}\label{lem_coll_cat_cupl}
For any $B\subset A \subset \R^N$ closed, we have
$$\cat(A,B) \geq \cupl(A,B)+1.$$
\end{lemma}

We are ready now to prove the main theorem.

\medskip

\claim Proof of Theorem \ref{claim:1.1}.
By construction of the neighborhood $\mc{X}_{\eps, \delta}$ and Corollary \ref{corol_passaggio_critici} (recall that $\rho_0 < r_3 \leq r_2'$ and that $J_{\eps}(u) < E_{m_0} + R(\hat{\delta}, u) \leq E_{m_0} + \hat{\delta} \leq l_0'$ for $u \in \mc{X}_{\eps, \delta}$), we have
$$
\Big\{ u \in (\mathcal{X}_{\eps, \delta})^{E_{m_0}+\hat{\delta}}_{E_{m_0}-\hat{\delta}} \mid J'_{\eps}(u)=0 \Big\} \subset \left\{ u \in H^s(\R^N) \mid I'_{\eps}(u)=0 \right\}.
$$
Thus we obtain
\begin{eqnarray*}
\lefteqn{ \# \{ u \textnormal{ solutions of \eqref{eq_princ_epsx}}\} \geq \# \Big\{ u \in (\mathcal{X}_{\eps, \delta})^{E_{m_0}+\hat{\delta}}_{E_{m_0}-\hat{\delta}} \mid J'_{\eps}(u)=0 \Big\} } \\
&\stackrel{(i)}\geq& \cat\Big (\mathcal{X}_{\eps, \delta}^{E_{m_0}+\hat{\delta}},\, \mathcal{X}_{\eps, \delta}^{E_{m_0}-\hat{\delta}}\Big) 
\stackrel{(ii)}\geq \cupl\Big (\mathcal{X}_{\eps, \delta}^{E_{m_0}+\hat{\delta}},\, \mathcal{X}_{\eps, \delta}^{E_{m_0}-\hat{\delta}}\Big) +1 \\
&\stackrel{(iii)}\geq& \cupl(K) +1
\end{eqnarray*}
that is the claim, up to the proof of (i)--(iii). Indeed, (i) is obtained classically from the Deformation Lemma \ref{lem_def_lem} (see e.g. \cite[Proposition 3.2]{Szu} or \cite[Theorem 6.1]{FLRW}). 
Inequality (ii) is given by the algebraic-topological Lemma \ref{lem_coll_cat_cupl}. 
Point (iii) is instead due to the existence of the homotopy gained in Proposition \ref{prop_esist_homot} and properties of the cup-length: indeed,
by \eqref{eq_propr_cupl} in Lemma \ref{lem_propr_cupl} (a), we have
$$\cupl \Big (\mathcal{X}_{\eps, \delta}^{E_{m_0}+\hat{\delta}},\, \mathcal{X}_{\eps, \delta}^{E_{m_0}-\hat{\delta}}\Big) \geq \cupl (\Psi_{\eps} \circ \Phi_{\eps});$$
moreover, since $\Psi_{\eps} \circ \Phi_{\eps}$ is homotopic to the immersion $j$ thanks to Proposition \ref{prop_esist_homot}, we have by Lemma \ref{lem_propr_cupl} (b)
$$ \cupl (\Psi_{\eps} \circ \Phi_{\eps}) = \cupl(j),$$
which leads to the conclusion thanks to Lemma \ref{lem_cupl_jK}. See Remark \ref{rem_regolar} for the proof of regularity.
\QED

\subsection{Concentration in $K$} \label{sec_concentration}

We prove now the polynomial decay and the concentration of the found solutions in $K$. 
First, we recall some definitions and properties of the fractional De Giorgi class introduced in \cite{Coz}, to which we refer for a complete introduction on the topic; we focus only on the linear case, where the exponent of the space is given by $2$.

Set first
$$\Tail(u; x_0, R):=(1-s) R^{2s} \int_{\R^N \setminus B_R(x_0)} \frac{|u(x)|}{|x-x_0|^{N+2s}}dx.$$
the \emph{tail function} of $u \in H^s(\R^N)$, centered in $x_0 \in \R^N$ with radius $R>0$, introduced in \cite{DKP1, DKP2}. From \cite[Paragraph 6.1]{Coz} we have the following definition.

\begin{definition}
Let $A \subset \R^N$ be open, $\zeta \geq 0$, $H\geq 1$, $k_0 \in \R$, $\mu \in (0, 2s/N]$, $\lambda \geq 0$ and $R_0 \in (0, +\infty]$. We say that $u$ belongs to the \emph{fractional De Giorgi class} $DG_+^{s,2}(A, \zeta, H, k_0, \mu , \lambda, R_0)$ if and only if
\begin{eqnarray*}
\lefteqn{[(u-k)_+]^{2}_{B_r(x_0)}+ \int_{B_r(x_0)} (u(x)-k)_+ \Bigg( \int_{B_{2R_0}(x)} \frac{(u(y)-k)_{-}}{|x-y|^{N+2s}} dy\Bigg) dx } \\
&\leq & \frac{H}{1-s} \bigg(\Big(R^{\lambda} \zeta^{2} + \frac{|k|^{2}}{R^{N \mu }}\Big) \big| \supp((u-k)_+) \cap B_R(x_0) \big |^{1- \frac{2s}{N} + \mu } +\\
&&+ \frac{R^{2(1-s)}}{(R-r)^{2}} \norm{(u-k)_+}_{L^{2}(B_R(x_0))}^{2} + \\
&& + \frac{R^{N }}{(R-r)^{N+2s}} \norm{(u-k)_+}_{L^1(B_R(x_0))} \Tail((u-k)_+; x_0, r) \bigg)
\end{eqnarray*}
for any $x_0 \in A$, $0<r<R< \min\{R_0, d(x_0, \partial A)\}$ and $k\geq k_0$. Here $(u-k)_{\pm}$ denote the positive and negative parts of the function $u-k$.
\end{definition}

We see now how this class of functions is related to the PDE setting. By a careful analysis of the proof of \cite[Proposition 8.5]{Coz} we obtain the following result.

\begin{theorem}\label{thm_degiorgi_class}
Let $N\geq 2$ and let $u\in H^s(\R^N)$ be a weak subsolution of
\begin{equation}\label{eq_eq_cozzi}
(-\Delta)^s u \leq g(x,u), \quad x \in \R^N
\end{equation}
where $g: \R^N \times \R \to \R$ satisfies, for a.e. $x \in \R^N$ and every $t \in \R$,
$$\abs{g(x,t)} \leq d_1 + d_2 |t|^{q-1}$$
for some $q \in (2, 2^*_s)$. Then there exist $\alpha=\alpha(N,s,q)>0$, $C=C(N,s,q,d_2)>0$ and $H=H(N, s, q, d_2)\geq 1$ such that, for each $x_0\in \R^N$ and each $R_0$ verifying
$$0<R_0 \leq C(N,s,q,d_2) \min \left \{1, \norm{u}_{L^{2^*_s}(\R^N)}^{-\alpha(N,s,q)}\right\},$$
it results that
$$u \in \DG_+^{s,2} \Big( B_{R_0}(x_0), d_1, H, 0, 1-\frac{q}{2^*_s}, 2s, R_0\Big).$$
\end{theorem}

As shown in \cite[Proposition 6.1 and Theorem 8.2]{Coz}, the belonging to a De Giorgi class implies useful $L^{\infty}_{loc}$ and $C^{0,\sigma}_{loc}$ estimates, which will be implemented in the following proof.

\medskip

\claim Proof of Theorem \ref{mean-concen}. 
For $\eps$ sufficiently small, let $u_{\eps}$ be one of the $\cupl(K)+1$ critical points of $J_{\eps}$ built in Theorem \ref{mean-concen}, which by Corollary \ref{corol_passaggio_critici} is also a solution of \eqref{eq_princ_epsx}, positive by (f2). In particular, since it satisfies the assumptions of Lemma \ref{lemma_stima_basso}, looking at the proof (see \eqref{eq_conv_uu} and \eqref{eq_conv_p_eps}) we obtain that
$$\norm{u_{\eps} - U(\cdot - p_{\eps})}_{H^s(\R^N)}\to 0$$
with $U \in S_{V(p_0)}$, $p_{\eps}\in \R^N$ and
$$\eps p_{\eps} \to p_0 \in \Omega[0, \nu_1).$$

\textbf{Step 1.}
Notice that we have found these solutions by fixing $\nu_0$, $l_0$ and $l_0'$. Let them move, throughout three sequences $\nu_0^n \searrow 0$, $l_0^n \searrow E_{m_0}$, and $(l_0')^n \searrow E_{m_0}$, 
and find the corresponding (sufficiently small) $\eps_n >0$ such that $\cupl(K)+1$ solutions exist; let $u_{\eps_n}$ be one of those and $p_{\eps_n}$ as before. It is not reductive to assume $\eps_n \to 0$ as $n \to +\infty$; by a diagonalization-like argument we obtain
\begin{equation} \label{eq_dim_conv_u_x}
u_{\eps_n}(\cdot+p_{\eps_n})\to U \; \; \textnormal{ in $H^s(\R^N)$}, \; \textnormal{ for some $U$ least energy solution of \eqref{eq_least_energy_m0}},
\end{equation}
$$\eps_n p_{\eps_n} \to p_0 \in K,$$
as $n \to +\infty$.

\textbf{Step 2.} \label{pag_dim_concent}
From now on we write $\eps\equiv \eps_n$ to avoid cumbersome notation. 
By $I_{\eps}'(u_{\eps})=0$ we obtain
\begin{equation}\label{eq_dim_perconfr}
(-\Delta)^s u_{\eps} + V(\eps x) u_{\eps} = f(u_{\eps}) , \quad x \in \R^N,
\end{equation}
thus (recall that $u_{\eps}$ is positive), by choosing $\beta < \underline{V}$ in \eqref{eq_prop_f}, 
$$(-\Delta)^s u_{\eps} \leq -\underline{V} u_{\eps} + f(u_{\eps})\leq (\beta-\underline{V})u_{\eps} + C_{\beta} u_{\eps}^p \leq C_{\beta} u_{\eps}^p, \quad x \in \R^N.$$
Therefore by Theorem \ref{thm_degiorgi_class} we have, choosing $q=p$, $d_1=0$ and $d_2=C_{\beta}$,
$$u_{\eps} \in \DG_+^{s,2} \Big( B_{R_0}(x_0), 0, H, 0, 1-\frac{p}{2^*_s}, 2s, R_0\Big),$$
with $H=H(N, s, p, \beta)$ 
and $R_0$ depending on $N, s, q, C_{\beta}$ and a uniform upper bound of the $H^s$-norms of $u_{\eps }$.
 
We can thus use now \cite[Proposition 6.1]{Coz}: observing that $d(x_0, \partial B_{R_0}(x_0))= R_0$, and that $\mu = 1-\frac{p}{2^*_s}$, we obtain, for any $\omega\in (0,1]$ and $R\in (0, \frac{R_0}{2})$, 
$$\sup_{B_R(x_0)} u_{\eps} \leq \frac{C}{(N-2s)^{\frac{1}{2\mu }}} \frac{1}{\omega^{\frac{1}{2\mu }}} \frac{1}{(2R)^{N/2}} \norm{u_{\eps}}_{L^2(B_{2R}(x_0))} + \omega \Tail(u_{\eps};x_0, R)$$
that is, rewriting the constant $C=C(N, s, p, \beta)$,
$$\sup_{B_R(x_0)} u_{\eps} \leq C\frac{1}{\omega^{\frac{1}{2\mu }}} \frac{1}{R^{N/2}} \norm{u_{\eps}}_{L^2(B_{2R}(x_0))} + \omega \Tail(u_{\eps};x_0, R).$$

\textbf{Step 3.}
We have
\begin{eqnarray*}
\norm{u_{\eps}}_{L^{\infty}(\R^N)} &=& \sup_{x_0 \in \R^N} \sup_{B_R(x_0)} u_{\eps} \\
&\leq& \sup_{x_0 \in \R^N} \left(C\frac{1}{\omega^{\frac{1}{2\mu }}} \frac{1}{R^{N/2}} \norm{u_{\eps}}_{L^2(B_{2R}(x_0))} + \omega \Tail(u_{\eps};x_0, R)\right).
\end{eqnarray*}
Observe that, by H\"older inequality,
\begin{eqnarray*}
\Tail(u_{\eps};x_0, R) &\leq& (1-s) R^{2s} \norm{u_{\eps}}_{L^2(\R^N \setminus B_R(x_0))} \norm{\tfrac{1}{|x-x_0|^{2N+4s}}}_{L^2(\R^N \setminus B_R(x_0))} \\
&\leq& \frac{C}{R^{N/2}}\norm{u_{\eps}}_{L^2(\R^N)}.
\end{eqnarray*}
Thus
\begin{eqnarray*}
\norm{u_{\eps}}_{L^{\infty}(\R^N)} &\leq& \frac{C}{R^{N/2}}\sup_{x_0 \in \R^N} \left(\omega^{-\frac{1}{2\mu }} \norm{u_{\eps}}_{L^2(B_{2R}(x_0))} + \omega \norm{u_{\eps}}_{L^2(\R^N \setminus B_R(x_0))}\right)\\
&\leq& \frac{C}{R^{N/2}}\left(\omega^{-\frac{1}{2\mu }} + \omega\right) \norm{u_{\eps}}_{L^2(\R^N)}
\end{eqnarray*}
which is uniformly bounded by the properties on $u_{\eps}$. Hence $u_{\eps}$ are uniformly bounded in $L^{\infty}(\R^N)$.

In addition, by the estimates on $V$, $f$ and $u_{\eps}$, we have 
$$g_{\eps}(x):= -V(\eps x) u_{\eps}(x) + f(u_{\eps}(x)) \in L^{\infty}(\R^N)$$
with bound uniform in $\eps$; since
$$(-\Delta)^s u_{\eps} = g_{\eps}(x), \quad x \in \R^N,$$
by \cite[Theorem 8.2]{Coz} there exists $\sigma \in (0,1)$, not depending on $u_{\eps}$, and $C=C(N,s)$, such that, for each $R>1$ and $x_0 \in \R^N$, 
\begin{eqnarray*}
\lefteqn{[u_{\eps}]_{C^{0,\sigma}(B_R(x_0))}}\\
&\leq& \frac{C}{R^{\sigma}} \left( \norm{u_{\eps}}_{L^{\infty}(B_{4R}(x_0))} + \Tail(u_{\eps}; x_0, 4R) + R^{2s} \norm{g_{\eps}}_{L^{\infty}(B_{8R}(x_0))}\right)
\leq C'.\label{eq_stima_uniforme_C0alpha}
\end{eqnarray*}
We highlight that, since the constant is uniform in $R>1$, we obtain $u_{\eps} \in C^{0, \sigma}(\R^N)$.

\textbf{Step 4.}
By the local uniform estimate on $u_{\eps}$ we could gain $\norm{u_{\eps}}_{L^{\infty}(\R^N\setminus \left(\Omega_{h_0}/\eps)_{2R_0}\right)}\to 0$, but this lack of uniformity on the domain can be improved. 
Thus we exploit the tightness of $\tilde{u}_{\eps}$ to reach the claim, where
$$\tilde{u}_{\eps}:=u_{\eps}(\cdot + p_{\eps}).$$
Indeed, by Step 2, and \eqref{eq_dim_conv_u_x} we have
$$\parag{ &(-\Delta)^s \tilde{u}_{\eps} +V(\eps x+ \eps p_{\eps}) \tilde{u}_{\eps} = f(\tilde{u}_{\eps}), \quad x \in \R^N,& \\ &\norm{\tilde{u}_{\eps}}_{\infty}\leq C,& \\ &\tilde{u}_{\eps} \to U \quad \textnormal{ in $H^s(\R^N)$ as $\eps \to 0$}, \quad \textnormal{ $U$ least energy solution of \eqref{eq_least_energy_m0}}.&}$$
In particular, it is standard to show that $f(\tilde{u}_{\eps}) \to f(U)$ in $L^2(\R^N)$, $\norm{f(\tilde{u}_{\eps})}_{\infty}\leq C$ and $U, f(U) \in L^{\infty}(\R^N)$. By interpolation we thus obtain
$$\chi_{\eps}:=\tilde{u}_{\eps} + f(\tilde{u}_{\eps})\to \chi:=U+f(U) \quad \textnormal{ in $L^{q'}(\R^N)$}$$
for every $q'\in [2, +\infty)$, and $\norm{\chi_{\eps}}_{\infty} \leq C.$ 
Proceeding as in \cite[Lemma 2.6]{AlMi} we obtain
\begin{equation}\label{eq_conv_unif_0_veps}
\tilde{u}_{\eps}(x) \to 0 \; \; \textnormal{ as $|x|\to +\infty$}, \quad \textnormal{uniformly in $\eps$}.
\end{equation}
For the reader's convenience, we sketch here the proof. Indeed, being $\tilde{u}_{\eps}$ solution of
$$(-\Delta)^s \tilde{u}_{\eps} + \tilde{u}_{\eps} = \chi_{\eps} - V(\eps x + \eps p_{\eps}) \tilde{u}_{\eps}, \quad x \in \R^N,$$
we have the representation formula
$$\tilde{u}_{\eps}= \mc{K} * (\chi_{\eps} - V(\eps x + \eps p_{\eps}) \tilde{u}_{\eps})$$
where $\mc{K}$ is the Bessel Kernel; we recall that $\mc{K}$ is positive, it satisfies $\mc{K}(x) \leq \frac{C}{|x|^{N+2s}}$ for $|x| \geq 1$ and $\mc{K} \in L^q(\R^N)$ for $q \in [1, 1 + \tfrac{2s}{N-2s})$ (see \cite[page 1241 and Theorem 3.3]{FQT}). Let us fix $\eta>0$; since $V$, $\tilde{u}_{\eps}$ and $\mc{K}$ are positive, we have, for $x \in \R^N$,
\begin{eqnarray*}
\tilde{u}_{\eps}(x) &=& \int_{\R^N} \mc{K}(x-y) \big(\chi_{\eps}(y) - V(\eps x + \eps p_{\eps})\tilde{u}_{\eps}(y)\big)dy \\
&\leq& \int_{|x-y|\geq 1/\eta} \mc{K}(x-y) \chi_{\eps}(y)dy +\int_{|x-y|< 1/\eta} \mc{K}(x-y) \chi_{\eps}(y)dy.
\end{eqnarray*}
As regards the first piece
$$ \int_{|x-y|\geq 1/\eta} \mc{K}(x-y) \chi_{\eps}(y)dy \leq \norm{\chi_{\eps}}_{\infty} \int_{|x-y|\geq 1/\eta} \frac{C}{|x-y|^{N+2s}} dy \leq C \eta^{2s}$$
while for the second piece, fixed a whatever $q \in (1, 1 + \tfrac{2s}{N-2s})$ and its conjugate exponent $q'> \frac{N}{2s}$, we have by H\"older inequality
\begin{eqnarray*}
\int_{|x-y|< 1/\eta} \mc{K}(x-y) \chi_{\eps}(y)dy &\leq& \norm{\mc{K}}_q \norm{\chi_{\eps}}_{L^{q'}(B_{1/\eta}(x))}\\
& \leq& \norm{\mc{K}}_q\left( \norm{\chi_{\eps}-\chi}_{q'}+ \norm{\chi}_{L^{q'}(B_{1/\eta}(x))}\right)
\end{eqnarray*}
where the first norm can be made small for $\eps<\eps_0=\eps_0(\eta)$, while the second 
 for $|x| \gg 0$ (uniformly in $\eps$). On the other hand, for $\eps\geq \eps_0$ (and thus for a finite number of elements, since we recall we are working with $\eps \equiv \eps_n$ small) the quantity $ \norm{\chi_{\eps}}_{L^{q'}(B_{1/\eta}(x))}$ can be made small for $|x|\gg 0$, uniformly in $\eps$. Joining the pieces, we have \eqref{eq_conv_unif_0_veps}.

\textbf{Step 5.}
Let now $y_{\eps}\in \R^N$ be a maximum point for $u_{\eps}$, which exists by the boundedness of $u_{\eps}$ and its continuity (see \eqref{eq_stima_uniforme_C0alpha}). 
Therefore $z_{\eps}:=y_{\eps}-p_{\eps}$ 
 is a maximum point for $\tilde{u}_{\eps}$. In particular
$$\tilde{u}_{\eps}(z_{\eps})=\max_{\R^N} \tilde{u}_{\eps}=\norm{\tilde{u}_{\eps}}_{\infty}\not \to 0 \quad \textnormal{ as $\eps \to 0$}$$
since on the contrary we would have $\tilde{u}_{\eps}\to 0$ almost everywhere, which is in contradiction with the fact that $\tilde{u}_{\eps}\to U\nequiv 0$ almost everywhere (up to a subsequence). As a consequence, thanks to \eqref{eq_conv_unif_0_veps}, we have that $z_{\eps}$ is bounded (up to a subsequence). That is, again up to a subsequence,
$$z_{\eps} \to \overline{p}$$
for some $\overline{p}\in \R^N$. In particular
$$\eps y_{\eps} = \eps z_{\eps} + \eps p_{\eps} \to p_0 \in K$$
and, by the fact that
$$U(\cdot + z_{\eps}) \to U(\cdot + \overline{p}) =: \overline{U} \quad \textnormal{ in $H^s(\R^N)$}$$
we have 
$u_{\eps}(\cdot + y_{\eps}) \to \overline{U}$ in $H^s(\R^N)$, $\overline{U}$ least energy solution of \eqref{eq_least_energy_m0}. 
We set
$$\overline{u}_{\eps}:= u_{\eps}(\cdot + y_{\eps}),$$
$$\overline{u}_{\eps}\to \overline{U} \; \; \textnormal{ in $H^s(\R^N)$}, \quad \textnormal{ $\overline{U}$ least energy solution of \eqref{eq_least_energy_m0}};$$
in addition, $\overline{u}_{\eps}$ is positive by \textnormal{(f2)}, and in the same way we obtained \eqref{eq_conv_unif_0_veps} we obtain also
\begin{equation}\label{eq_conv_unif_0_veps2}
\overline{u}_{\eps}(x) \to 0 \; \; \textnormal{ as $|x|\to +\infty$}, \quad \textnormal{uniformly in $\eps$}.
\end{equation}
Moreover, by exploiting the uniform estimates in $L^{\infty}(\R^N)$ and $C^{0, \sigma}_{loc}(\R^N)$ we obtain by Ascoli-Arzell\`{a} theorem also that $\overline{u}_{\eps}\to \overline{U}>0$ in $L^{\infty}_{loc}(\R^N)$, with $\overline{U}$ continuous; this easily implies, for every $r>0$, that
\begin{equation}\label{eq_stima_basso_min}
\min_{B_r} \overline{u}_{\eps} \geq \frac{1}{2} \min_{B_r} \overline{U} >0
\end{equation}
for $\eps$ small, depending on $\overline{U}$ and $r$.

\textbf{Step 6.}
By \eqref{eq_conv_unif_0_veps2} we have, for $R'$ large (uniform in $\eps$), that
$$\overline{u}_{\eps}(x)\leq \eta', \quad \textnormal{ for $|x|>R'$}$$
for every $\eps>0$, where $\eta'>0$ is preliminary fixed. 
As a consequence, by (f1.2), we gain
$$-\frac{1}{2} \overline{V} \overline{u}_{\eps}(x)\leq f(\overline{u}_{\eps}(x)) \leq \frac{1}{2}\underline{V} \overline{u}_{\eps}(x),\quad \textnormal{ for $|x|>R'$},$$
where $\overline{V}:=\norm{V}_{\infty}$. We obtain by \eqref{eq_dim_perconfr}
$$(-\Delta)^s \overline{u}_{\eps} + \tfrac{1}{2}\underline{V} \overline{u}_{\eps} \leq f(\overline{u}_{\eps}) - \tfrac{1}{2} \underline{V}\overline{u}_{\eps} \leq 0, \quad x \in \R^N \setminus B_{R'},$$
$$(-\Delta)^s \overline{u}_{\eps} + \tfrac{3}{2} \overline{V} \overline{u}_{\eps} \geq f(\overline{u}_{\eps}) + \tfrac{1}{2} \overline{V}\overline{u}_{\eps}\geq 0, \quad x \in \R^N \setminus B_{R'}.$$
Notice that we always intend differential inequalities in the weak sense, that is tested with functions in $H^s(\R^N)$ with supports contained (e.g.) in $\R^N \setminus B_{R'}$. 

In addition, by Lemma \ref{lem_esist_sol_part} we have that there exist two positive functions $\underline{W}'$, $\overline{W}'$ and three positive constants $R''$, $C'$ and $C''$ depending only on $V$, such that
$$ \parag{
& (-\Delta)^s \underline{W}' + \frac{3}{2}\overline{V} \, \underline{W}' = 0, \quad x \in \R^N \setminus B_{R''},& \\ 
&\frac{C'}{|x|^{N+2s}}< \underline{W}' (x), \quad \textnormal{ for $|x|>2R''$}.&}$$
and
$$ \parag{& (-\Delta)^s \overline{W}' + \frac{1}{2}\underline{V} \, \overline{W}' = 0, \quad x \in \R^N \setminus B_{R''},& \\ 
& \overline{W}'(x) < \frac{C''}{|x|^{N+2s}}, \quad \textnormal{ for $|x|>2R''$},&}$$
Set $R:=\max\{ R', 2R''\}$. Let $\underline{C}_1$ and $\overline{C}_1$ be some uniform lower and upper bounds for $\overline{u}_{\eps}$ on $B_R$, $\underline{C}_2:=\min_{B_R} \overline{W}'$ and $\overline{C}_2:= \max_{B_R} \underline{W}'$, all strictly positive. Define
$$\underline{W}:= \underline{C}_1 \overline{C}_2 ^{-1} \underline{W}', \quad \overline{W}:= \overline{C}_1 \underline{C}_2^{-1} \overline{W}'$$
so that
$$\underline{W}\leq \overline{u}_{\eps} \leq \overline{W}, \quad \textnormal{ for $|x|\leq R$}.$$
Through a Comparison Principle (see Lemma \ref{lem_comp_prin}), and redefining $C'$ and $C''$, we obtain
$$ \frac{C'}{|x|^{N+2s}} <\underline{W}(x) \leq\overline{u}_{\eps}(x) \leq \overline{W}(x) < \frac{C''}{|x|^{N+2s}}, \quad \textnormal{ for $|x|>R$}.$$
By the uniform boundedness of $\overline{u}_{\eps}$ and \eqref{eq_stima_basso_min} we also obtain
$$\frac{C'}{1+|x|^{N+2s}}< \overline{u}_{\eps}(x) < \frac{C''}{1+|x|^{N+2s}}, \quad \textnormal{ for $x \in \R^N $}.$$
Recalling the definition of $\overline{u}_{\eps}$, we have finally obtained a sequence of solutions such that
$$\parag{&u_{\eps_n}(y_{\eps_n})=\max_{\R^N} u_{\eps_n}, &\\& 
d(\eps_n y_{\eps_n}, K)\to 0, &\\& 
\frac{C'}{1+|x-y_{\eps_n}|^{N+2s}} \leq u_{\eps_n}(x) \leq \frac{C''}{1+|x-y_{\eps_n}|^{N+2s}}, \quad \textnormal{ for $x \in \R^N$}, &\\& 
\norm{u_{\eps_n}(\cdot + y_{\eps_n}) - \overline{U}}_{H^s(\R^N)} \to 0, \quad \textnormal{ for some $\overline{U}$ least energy solution of \eqref{eq_least_energy_m0} },}$$
where the limits are given by $n\to +\infty$. Furthermore, by the uniform estimates in $L^{\infty}(\R^N)$ and the local uniform estimates in $C^{0, \sigma}_{loc}(\R^N)$ of $u_{\eps_n}$, together with the locally-compact version of Ascoli-Arzel\`{a} theorem, we have that the last convergence is indeed uniform on compacts. 
Thus, recalled that $v_{\eps_n}= u_{\eps_n}(\cdot/\eps_n)$ are solutions of the original problem \eqref{eq:1.1}, defined $x_{\eps_n}:=\eps_n y_{\eps_n}$ we obtain, as $n \to +\infty$,
 $$\parag{&v_{\eps_n}(x_{\eps_n})=\max_{\R^N} v_{\eps_n}, &\\& 
 d(x_{\eps_n}, K)\to 0, &\\& 
 \frac{C'}{1+|\frac{x-x_{\eps_n}}{\eps_n}|^{N+2s}}\leq v_{\eps_n}(x) \leq \frac{C''}{1+|\frac{x-x_{\eps_n}}{\eps_n}|^{N+2s}}, \quad \textnormal{ for $x \in \R^N$}, &\\& 
 \norm{v_{\eps_n}(\eps_n \cdot + x_{\eps_n}) - \overline{U}}_{X} \to 0, \quad X= H^s(\R^N) \textnormal{ and } X=L^{\infty}_{loc}(\R^N), }$$
for some $\overline{U}$ least energy solution of \eqref{eq_least_energy_m0}. This concludes the proof.
\QED

\begin{remark}\label{rem_regolar}
We observe that Steps 2 and 3 apply to a whatever family of equations $(u_{\eps})_{\eps}$, that is why the regularity statement in Theorem \ref{claim:1.1} holds true.
\end{remark}

\appendix

\section{Polynomial decay of $\widehat{S}$} \label{sez_poly_dec}

We recall some lemmas useful to prove polynomial decay estimates. We start with a comparison principle: the result is a straightforward adaptation of \cite[Lemma 6]{SV}, and we prove it here for the sake of completeness. See also \cite{Jar}.

\begin{lemma}[Maximum Principle]\label{lem_comp_prin}
Let $\Sigma \subset \R^N$, possibly unbounded, and let $u\in H^s(\R^N)$ be a weak subsolution of
$$(-\Delta)^s u + a u \leq 0, \quad x \in \R^N \setminus \Sigma$$
with $a>0$, in the sense that
$$\int_{\R^N} (-\Delta)^{s/2} u \,(-\Delta)^{s/2} v \, dx+ a \int_{\R^N} u v \, dx\leq 0$$
for every positive $v \in H^s(\R^N)$ with $\supp(v) \subset \R^N \setminus \Sigma$. 
Assume moreover that
$$u\leq 0, \quad \textit{ for a.e. $x \in \Sigma$}.$$
Then
\begin{equation}\label{eq_dis_comp_prin}
u\leq 0, \quad \textit{ for a.e. $x \in \R^N$}.
\end{equation}
\end{lemma}

\claim Proof.
By the assumption we have $u^+=0$ on $\Sigma$, thus $u^+\in H^s(\R^N)$ is a suitable test function and we obtain, using $u=u^+-u^-$ and $u^+ u^-\equiv 0$,
\begin{eqnarray*}
0 &\geq& \int_{\R^N} |(-\Delta)^{s/2} u^+|^2 \, dx + a \int_{\R^N} |u^+|^2 \, dx - \int_{\R^N} (-\Delta)^{s/2} u^- (-\Delta)^{s/2} u^+ \, dx \\
&=& \norm{(-\Delta)^{s/2}u^+}_2^2 + a\norm{u^+}_2^2+C\int_{\R^{2N}} \frac{u^-(x)u^+(y)+u^-(y)u^+(x)}{|x-y|^{N+2s}} \, dx \, dy \\
&\geq& \norm{(-\Delta)^{s/2}u^+}_2^2 + a\norm{u^+}_2^2
\end{eqnarray*}
which implies $u^+=0$.
\QED

\begin{remark}
We point out that if $u$ is assumed continuous, then \eqref{eq_dis_comp_prin} is actually pointwise. Moreover, the constant $a>0$ may be substituted by a more general $a(x)>0$ which gives sense to the integrals. 
\end{remark}

The following result can be obtained similarly to \cite[Lemmas 4.2 and 4.3]{FQT}. 

\begin{lemma}[Comparison function]\label{lem_esist_sol_part}
Let $b>0$. Then there exists a strictly positive continuous function $W_b\in H^s(\R^N)$ such that, for some positive constants $C'_b, C_b''$, it verifies 
$$(-\Delta)^s W_b + \frac{b}{2} W _b = 0,\quad x \in \R^N \setminus B_{r_b}$$
pointwise, with $r_b:= \left(\frac{2}{b}\right)^{1/2s}$, and
\begin{equation}\label{eq_stima_fun_confr}
\frac{C'_b}{|x|^{N+2s}}<W_b(x)< \frac{C_b''}{|x|^{N+2s}}, \quad \textit{ for $|x|>2 r_b$}.
\end{equation}
The constants $r_b, C'_b, C_b''$ remain bounded by letting $b$ vary in a compact set far from zero.
\end{lemma}

\claim Proof.
Let $B_{1/2}\prec \varphi \prec B_1$, and define $\tilde{W}:= \mc{K} * \varphi$, where $\mc{K}$ is the Bessel potential. Arguing as in \cite{FQT} (see also \cite[Theorem 1.3]{BKS}) we obtain
$$(-\Delta)^s \tilde{W} + \tilde{W} = \varphi, \quad x \in \R^N$$
and
$$\frac{C'}{|x|^{N+2s}}<\tilde{W}(x)\leq \frac{C''}{|x|^{N+2s}}.$$
By scaling $W:=\tilde{W}(r_b \cdot)$ we reach the claim.
\QED

\begin{proposition}[Polynomial decay]\label{prop_polyn_dec}
Assume \textnormal{(f1)--(f3)}. Let $a>0$ and let $U$ be a weak solution of
$$(-\Delta)^s U + a U = f(U), \quad x \in \R^N.$$
Then there exist positive constants $C_a', C''_a$ such that
$$\frac{ C_a'}{1+|x|^{N+2s}} \leq U(x)\leq \frac{ C_a''}{1+|x|^{N+2s}},\quad \textit{ for $x \in \R^N$}.$$
These constants can be chosen uniform for $U \in \widehat{S}$.
\end{proposition}

\claim Proof. 
The proof is almost the same of the one carried out in Theorem \ref{mean-concen}.

Indeed, as in Step 2 and Step 3, we obtain the uniform boundedness in $L^{\infty}(\R^N)$. We point out that the values $C_{\delta}$, $H$, $C$ and $R_0$ depend on $a \in [m_0, m_0+\nu_0]$, since they depend on $\delta$ and we must have $\delta < a$; on the other hand, it is sufficient to take $\delta < m_0$ to gain uniformity. 
The same can be said on the uniform boundedness in $C^{0, \sigma}(\R^N)$ and for the constants $R'', C', C''$ related to the comparison functions $\underline{W}', \overline{W}'$, thanks to Lemma \ref{lem_esist_sol_part}. 
As we will show, this allows us to gain that
\begin{equation}\label{eq_conv_unif_0}
\lim_{|x|\to +\infty} U(x)=0 \quad \textnormal{ uniformly for $U\in \widehat{S}$},
\end{equation}
which leads, as in Step 6 of the proof, to
$$\abs{f(U(x))} \leq \frac{1}{2}a U(x),\quad \textnormal{ for $|x|>R'$}$$
where $R'$ does not depend on $a\in [m_0, m_0+\nu_0]$. In addition, compactness of $\widehat{S}$ and a simple contradiction argument lead to $\min_{B_R} U \geq C>0$ uniformly for $U\in \widehat{S}$. If we prove \eqref{eq_conv_unif_0}, we conclude as in Step 6. 

Let us prove \eqref{eq_conv_unif_0}. By contradiction, there exist $(x_k)_k\subset \R^N$, $|x_k|\to +\infty$, $(U_k)_k \subset \widehat{S}$ and $\theta>0$ such that $U_k(x_k)> \theta>0$. Define
$$V_k:= U_k(\cdot+x_k).$$
Since both are bounded sequences in $H^s(\R^N)$, we have $U_k \wto U$ and $V_k \wto V$ in $H^s(\R^N)$; moreover, by the uniform $L^{\infty}(\R^N)$ and $C^{0,\sigma}_{loc}(\R^N)$ estimates and Ascoli-Arzel\`{a} theorem, we have also that the convergences are pointwise. In particular by
$$U_k(0)\geq U_k(x_k) > \theta, \quad V_k(0)=U_k(x_k)>\theta$$
we obtain
$$U(0)\geq \theta >0, \quad V(0) \geq \theta>0.$$
As a consequence, $U$ and $V$ are not trivial. Let now $(a_k)_k \subset \R$ be such that $U_k \in S_{a_k}$; up to a subsequence we have $a_k \to a\in [m_0, m_0+\nu_0]$. Observed that also $V_k$ are solutions of $L'_{a_k}(V_k)=0$, we obtain, as in Step 3 of Lemma \ref{lemma_stima_unif} (see also Step 7 of the proof of Lemma \ref{lemma_stima_basso}), that $U$ and $V$ are (nontrivial) solutions of $L_a'(U)=0$. Hence
$$E_{m_0} \leq E_a \leq L_a(U), \quad E_{m_0} \leq E_a\leq L_a(V).$$
By the Pohozaev identity we have the following chain of inequalities, once fixed $R>0$ and $k\gg 0$ such that $|x_k|\geq 2R$,
\begin{eqnarray*}
\lefteqn{l_0 \geq \liminf_{k\to +\infty} L_{a_k}(U_k) = \frac{s}{N} \liminf_{k\to +\infty} \int_{\R^N}\abs{ (-\Delta)^{s/2}U_k}^2 dx} \\
&\geq& \frac{s}{N} \liminf_{k\to +\infty} \Bigg( \int_{B_R} \abs{ (-\Delta)^{s/2}U_k}^2 dx + \int_{B_R} \abs{ (-\Delta)^{s/2}V_k}^2 dy \Bigg) \\
&\geq & \frac{s}{N} \Bigg( \int_{B_R} \abs{ (-\Delta)^{s/2}U}^2 dx + \int_{B_R} \abs{ (-\Delta)^{s/2}V}^2 dy\Bigg)
\end{eqnarray*}
where in the last passage we have used that $U_k \wto U$ in $H^s(\R^N)$, thus $(-\Delta)^{s/2} U_k \wto (-\Delta)^{s/2} U$ in $L^2(\R^N)$, hence (by restriction) in $L^2(B_R)$, and the weak lower semicontinuity of the norm. 
Thus, by choosing $R$ sufficiently large, we have 
\begin{eqnarray*}
l_0 &\geq& \frac{s}{N} \Bigg( \int_{\R^N} \abs{ (-\Delta)^{s/2}U}^2 dx+ \int_{\R^N} \abs{ (-\Delta)^{s/2}V}^2 dy\Bigg) - \eta\\
&=& L_a(U) + L_a(V) - \eta 
\geq 2E_{m_0}- \eta
\end{eqnarray*}
which leads to a contradiction if we choose $\eta \in (0, 2E_{m_0} - l_0)$, possible thanks to \eqref{eq_stima_l0}.
\QED

\bigskip
\bigskip

{\bf Acknowledgments.}
The authors are supported by MIUR-PRIN project ``Qualitative and quantitative aspects of nonlinear PDEs'' (2017JPCAPN\_005), and partially supported by GNAMPA-INdAM.

The authors would like also to thank Prof. K. Tanaka for the fruitful comments.

\end{document}